\documentclass[11pt]{article}

\usepackage[T1]{fontenc}
\usepackage{lmodern}
\usepackage{amsmath,amssymb,amsthm}
\usepackage{xcolor,tikz}
\usepackage[margin=1in]{geometry}
\usepackage[hidelinks]{hyperref}
\usepackage[nameinlink,noabbrev]{cleveref}
\usetikzlibrary{arrows.meta}
\usepackage{fontawesome5}

\usepackage[backend=biber,style=numeric]{biblatex}
\newtheorem{theorem}{Theorem}
\numberwithin{theorem}{section}

\newtheorem{lemma}[theorem]{Lemma}
\newtheorem{proposition}[theorem]{Proposition}
\newtheorem{corollary}[theorem]{Corollary}

\theoremstyle{definition}
\newtheorem{definition}[theorem]{Definition}

\crefname{theorem}{Theorem}{Theorems}
\Crefname{theorem}{Theorem}{Theorems}

\crefname{lemma}{Lemma}{Lemmas}
\Crefname{lemma}{Lemma}{Lemmas}

\crefname{proposition}{Proposition}{Propositions}
\Crefname{proposition}{Proposition}{Propositions}

\crefname{corollary}{Corollary}{Corollaries}
\Crefname{corollary}{Corollary}{Corollaries}

\crefname{definition}{Definition}{Definitions}
\Crefname{definition}{Definition}{Definitions}

\crefname{example}{Example}{Examples}
\Crefname{example}{Example}{Examples}

\crefname{remark}{Remark}{Remarks}
\Crefname{remark}{Remark}{Remarks}

\AddToHook{env/lemma/begin}{%
  \crefalias{theorem}{lemma}%
}

\AddToHook{env/proposition/begin}{%
  \crefalias{theorem}{proposition}%
}

\AddToHook{env/corollary/begin}{%
  \crefalias{theorem}{corollary}%
}

\AddToHook{env/definition/begin}{%
  \crefalias{theorem}{definition}%
}

\AddToHook{env/example/begin}{%
  \crefalias{theorem}{example}%
}

\AddToHook{env/remark/begin}{%
  \crefalias{theorem}{remark}%
}

\newcommand{\rmd}{\mathrm{d}}
\newcommand{\E}{\mathbb{E}}
\newcommand{\mBB}{\mathrm{mBB}}
\newcommand{\Mcal}{\mathcal{M}}
\newcommand{\Ecal}{\mathcal{E}}
\newcommand{\Pcal}{\mathcal{P}}
\newcommand{\Wcal}{\mathcal{W}}
\newcommand{\R}{\mathbb{R}}
\newcommand{\Rd}{\mathbb{R}^d}
\newcommand{\MCov}{\operatorname{MCov}}

\definecolor{phiblue}{RGB}{56,126,184}
\definecolor{psiorange}{RGB}{230,126,34}

\colorlet{philow}{phiblue!14}
\colorlet{phimid}{phiblue!28}
\colorlet{phihigh}{phiblue!42}
\colorlet{psilow}{psiorange!14}
\colorlet{psimid}{psiorange!28}
\colorlet{psihigh}{psiorange!42}
\colorlet{excludedcell}{gray!12}

\tikzset{
    plotframe/.style={black, line width=0.30pt},
    philabel/.style={font=\small, text=phiblue!75!black},
    psilabel/.style={font=\small, text=psiorange!75!black},
    pointlabel/.style={font=\scriptsize, inner sep=1pt},
    targetpoint/.style={
        circle, draw=black, line width=0.8pt,
        fill=none, inner sep=0pt, minimum size=5pt
    },
    arrowmass/.style={font=\scriptsize, inner sep=0.8pt},
    masssource/.style={arrowmass, left=2pt},
    firstarrow/.style={
        -{Latex[length=1.95mm, width=1.45mm]},
        draw=green!55!black, line width=0.95pt,
        line cap=round, shorten <=4.3pt, shorten >=4.2pt
    },
    shortfirstarrow/.style={
        -{Latex[length=1.95mm, width=1.45mm]},
        draw=green!55!black, line width=0.95pt,
        line cap=round, shorten <=3.0pt, shorten >=4.2pt
    },
    secondarrow/.style={
        -{Latex[length=1.95mm, width=1.45mm]},
        draw=violet!85!black, line width=0.95pt,
        line cap=round, shorten <=4.3pt, shorten >=4.2pt
    }
}

\newcommand{\circleat}[1]{%
    \node[targetpoint] at #1 {};
}
\newcommand{\crossat}[1]{%
    \begin{scope}[shift={#1}]
        \draw[black, line width=0.8pt]
            (-3.55pt,-3.55pt) -- (3.55pt,3.55pt)
            (-3.55pt,3.55pt) -- (3.55pt,-3.55pt);
    \end{scope}%
}

\newcommand{\targetcircles}{%
    \foreach \point in {
        (3,3),(3,4),(3,6),(4,6),(6,6),
        (6,8),(6,12),(8,12),(12,12),(12,16)
    } {\circleat{\point}}
}
\newcommand{\sourcecrosses}{%
    \crossat{({24/7},{30/7})}
    \crossat{({30/7},{48/7})}
    \crossat{({48/7},{60/7})}
    \crossat{({60/7},{96/7})}
}

\newcommand{\localframe}[1]{%
    \draw[plotframe] (2,2) rectangle (12,16);
    \foreach \xx/\lab in {3/{3#1},6/{6#1}} {
        \draw[black] (\xx,2) -- ++(0,-0.17)
            node[below, font=\scriptsize] {$\lab$};
    }
    \foreach \yy/\lab in {6/{6#1},12/{12#1}} {
        \draw[black] (2,\yy) -- ++(-0.17,0)
            node[left, font=\scriptsize] {$\lab$};
    }
    \node[below right, font=\small] at (12,2) {$s$};
    \node[above left, font=\small] at (2,16) {$t$};
}

\newcommand{\localregions}[1]{%
    \fill[#1]      (2,2)--(3,2)--(3,4)--(2,3)--cycle;
    \fill[phimid]  (3,2)--(6,2)--(6,8)--(4,6)--(3,4)--cycle;
    \fill[phihigh] (6,2)--(12,2)--(12,16)--(8,12)--(6,8)--cycle;
    \fill[psilow]  (2,3)--(3,4)--(4,6)--(2,6)--cycle;
    \fill[psimid]  (2,6)--(4,6)--(6,8)--(8,12)--(2,12)--cycle;
    \fill[psihigh] (2,12)--(8,12)--(12,16)--(2,16)--cycle;
}

\newcommand{\localregionlabels}[2]{%
    \node[philabel] at (4.55,2.95) {$\phi_{#1}$};
    \node[philabel] at (9.25,4.60) {$\phi_{#2}$};
    \node[psilabel] at (3.25,9.45) {$\psi_{#1}$};
    \node[psilabel] at (4.20,13.15) {$\psi_{#2}$};
}

\newcommand{\constructionframe}{%
    \draw[plotframe] (0,0) rectangle (10,14);
    \foreach \xx in {2,4,6,8} {
        \draw[black] (\xx,0) -- ++(0,-0.17)
            node[below, font=\scriptsize] {$\xx$};
    }
    \foreach \yy in {2,4,6,8,10,12} {
        \draw[black] (0,\yy) -- ++(-0.17,0)
            node[left, font=\scriptsize] {$\yy$};
    }
    \node[below right, font=\small] at (10,0) {$s$};
    \node[above left, font=\small] at (0,14) {$t$};
}

\newcommand{\constructionregions}{%
    \fill[excludedcell] (0,0) rectangle (10,14);
    \begin{scope}
        \clip (0,0) rectangle (10,14);
        \fill[phimid]
            (3,0)--(6,0)--(6,8)--(4,6)--(3,4)--cycle;
        \fill[phihigh]
            (6,0)--(12,0)--(12,16)--(8,12)--(6,8)--cycle;
        \fill[psilow]
            (0,3)--(2,3)--(3,4)--(4,6)--(0,6)--cycle;
        \fill[psimid]
            (0,6)--(4,6)--(6,8)--(8,12)--(0,12)--cycle;
        \fill[psihigh]
            (0,12)--(8,12)--(12,16)--(16,24)--(0,24)--cycle;
    \end{scope}
}

\newcommand{\constructionlabels}{%
    \localregionlabels{1}{2}%
    \node[psilabel] at (1.05,4.15) {$\psi_0$};
}

\newcommand{\originalmarks}{%
    \crossat{(3,4)}
    \circleat{(3,4)}
    \crossat{({24/7},{30/7})}
    \crossat{({30/7},{48/7})}
    \crossat{(5,1)}
    \crossat{(1,5)}
    \circleat{(3,3)}
    \circleat{(3,6)}
    \circleat{(4,6)}
    \circleat{(6,8)}
    \crossat{({48/7},{60/7})}
    \crossat{({60/7},{96/7})}
    \crossat{(10,2)}
    \crossat{(2,10)}
    \circleat{(6,6)}
    \circleat{(6,12)}
    \circleat{(8,12)}
}

\newcommand{\spreadhalo}[5]{%
    \begin{scope}
        \clip #4;
        \fill[#1, opacity=.18] (#2,#3) circle[radius={.38*(#5)}];
        \fill[#1, opacity=.13] (#2,#3) circle[radius={.78*(#5)}];
        \fill[#1, opacity=.09] (#2,#3) circle[radius={1.18*(#5)}];
        \fill[#1, opacity=.055] (#2,#3) circle[radius={1.58*(#5)}];
        \fill[#1, opacity=.028] (#2,#3) circle[radius={1.98*(#5)}];
        \fill[#1, opacity=.012] (#2,#3) circle[radius={2.35*(#5)}];
    \end{scope}
}

\newcommand{\addedhalos}[1]{%
    \spreadhalo{#1}{5}{1}{(3,0)--(6,0)--(6,8)--(4,6)--(3,4)--cycle}{1.5}
    \spreadhalo{#1}{10}{2}{(6,0)--(12,0)--(12,16)--(8,12)--(6,8)--cycle}{2*1.5}
    \spreadhalo{#1}{2}{10}{(0,6)--(4,6)--(6,8)--(8,12)--(0,12)--cycle}{2*1.5}
    \spreadhalo{#1}{1}{5}{(0,3)--(2,3)--(3,4)--(4,6)--(0,6)--cycle}{1.5}
    \spreadhalo{#1}{4}{20}{(0,12)--(8,12)--(12,16)--(16,24)--(0,24)--cycle}{4*1.5}
}

\newcommand{\atomlabels}{%
    \node[font=\scriptsize, inner sep=0pt, outer sep=0pt,
          above left=2.0pt] at ({30/7},{48/7}) {$X_2$};
    \node[font=\scriptsize, inner sep=0pt, outer sep=0pt,
          below right=2.0pt] at ({48/7},{60/7}) {$X_3$};
    \node[font=\scriptsize, inner sep=0pt, outer sep=0pt,
          below right=2.0pt] at (3,3) {$S_0$};
    \node[font=\scriptsize, inner sep=0pt, outer sep=0pt,
          above left=2.0pt] at (3,6) {$R_1$};
}

\newcommand{\Tpt}[2]{%
    ({(#1)/(1+(#1)+(#2))},{(#2)/(1+(#1)+(#2))})%
}

\newcommand{\Tcrossat}[2]{%
    \pgfmathsetmacro{\tmarkx}{(#1)/(1+(#1)+(#2))}%
    \pgfmathsetmacro{\tmarky}{(#2)/(1+(#1)+(#2))}%
    \begin{scope}[shift={(\tmarkx,\tmarky)}, scale=0.53]
        \crossat{(0,0)}
    \end{scope}%
}

\newcommand{\Tcircleat}[2]{%
    \pgfmathsetmacro{\tmarkx}{(#1)/(1+(#1)+(#2))}%
    \pgfmathsetmacro{\tmarky}{(#2)/(1+(#1)+(#2))}%
    \circleat{(\tmarkx,\tmarky)}%
}

\newcommand{\setcellshade}[1]{%
    \pgfmathtruncatemacro{\cellshade}{
        ifthenelse(#1<=2,14*(#1+1),42+6*(#1-2))
    }%
}

\newcommand{\compactphicell}[1]{%
    \pgfmathsetmacro{\ascl}{pow(2,#1-1)}%
    \setcellshade{#1}%
    \fill[phiblue!\cellshade]
        \Tpt{3*\ascl}{0}--\Tpt{6*\ascl}{0}--
        \Tpt{6*\ascl}{8*\ascl}--\Tpt{4*\ascl}{6*\ascl}--
        \Tpt{3*\ascl}{4*\ascl}--cycle;
}

\newcommand{\compactpsicell}[1]{%
    \pgfmathsetmacro{\bscl}{pow(2,#1)}%
    \setcellshade{#1}%
    \fill[psiorange!\cellshade]
        \Tpt{0}{3*\bscl}--\Tpt{2*\bscl}{3*\bscl}--
        \Tpt{3*\bscl}{4*\bscl}--\Tpt{4*\bscl}{6*\bscl}--
        \Tpt{0}{6*\bscl}--cycle;
}

\title{Bass martingales need not stretch}
\author{M.~Beiglböck\thanks{\textit{Funding.} MB gratefully acknowledges financial support from FWF through P34743, P35197, and FW506064 and from OeNB through AB1898311.}, M.~Hasenbichler, G.~Pammer}
\date{}

\begin{document}

\maketitle

\begin{abstract}
    In classical optimal transport, the equivalence of optimality, convex gradient structure, and cyclical monotonicity is a central paradigm. 
    Known counterexamples to this equivalence rely on carefully constructed costs.
    In contrast, we show that this equivalence can fail in martingale transport in dimension~${d \geq 2}$: already for the canonical martingale Benamou--Brenier problem and irreducible marginals we have the strict inclusions
    \[
        \text{optimality}
        \;\subsetneq \;
        \text{(Bass) convex gradient structure}
        \;\subsetneq\;
        \text{martingale cyclical monotonicity}\, .
    \]
    However, under the stronger assumption of strict irreducibility, we establish that these three notions coincide.
    Understanding when the gradient structure guarantees optimality is also of practical relevance, since algorithms based on the Bass gradient representation are used in financial model calibration.
    
    All our counterexamples are based on the following phenomenon: for  $d\ge2$, there exist two one-step martingales $X,Y$ with the same initial and terminal laws, and a convex function $\psi$ 
    such that
    \[
        0=\mathbb E\big[\mathbb E[\psi(Y_1)\mid Y_0]-\psi(Y_0)\big]
        <
        \mathbb E\big[\mathbb E[\psi(X_1)\mid X_0]-\psi(X_0)\big]
        <\infty.
    \]
\end{abstract}

\noindent{\small\textbf{Keywords:} Martingale optimal transport; Bass martingales; stretched Brownian motion; cyclical monotonicity; convex duality; instability.\par}

\section{Introduction}\label{sec:introduction}

Throughout, unless stated otherwise, $\mu$ and $\nu$ are probability measures on $\R^d$ with finite second moments in convex order, $\mu\le_{\mathrm{cx}}\nu$.
By Strassen's theorem~\cite{St65}, the set $\Mcal(\mu,\nu)$ of martingale couplings between $\mu$ and $\nu$ is nonempty.
In contrast to classical transport, the marginals $\mu,\nu$ may split into components which cannot communicate through martingale couplings.
To exclude this, it is natural to assume that $(\mu,\nu)$ is \emph{irreducible}: for all Borel sets $A,B\subseteq\R^d$ with $\mu(A)>0$ and $\nu(B)>0$, there exists $\pi\in\Mcal(\mu,\nu)$ such that $\pi(A\times B)>0$; see, e.g.,~\cite[Definition~1.2]{BaBeScTs26}.

Stretched Brownian motion (sBM) was introduced in \cite{BaBeHuKa20,Lo18} and provides a martingale analogue of McCann's displacement interpolation~\cite{Mc97}.
It is defined as the unique optimizer of the martingale Benamou--Brenier problem
\begin{equation}\label{eq:mBB}\tag{mBB}
    \mBB(\mu,\nu)
    = \sup_{\substack{\rmd X_t = \sigma_t\,\rmd B_t, \\
                      X_0\sim\mu,\,X_1\sim\nu}}
      \E\Bigl[\int_0^1 \operatorname{tr}(\sigma_t)\,\rmd t\Bigr].
\end{equation}
Equivalently, sBM minimizes the adapted $2$-Wasserstein distance to Brownian motion, among martingales with initial law $\mu$ and terminal law $\nu$~\cite{BaBeHuKa20}.
The optimizer is unique in law and is a strong Markov martingale.
It yields a time-consistent martingale interpolation between $\mu$ and $\nu$: after rescaling time, its restriction to any subinterval is again sBM between the corresponding marginals~\cite{BaBeHuKa20}.
These properties motivate viewing sBM as a canonical martingale transport between measures in convex order, analogous to the Brenier coupling in classical transport.

A fundamental structural property of sBM is that it admits a representation as a Bass martingale~\cite[Theorem~1.3]{BaBeScTs26}.
A \emph{Bass martingale} is a process of the form
\begin{equation}\label{eq:bass}\tag{Bass}
    M_t=\E\bigl[\nabla\psi(B_1)\mid B_t\bigr],
    \qquad 0\le t\le1,
\end{equation}
where $B$ is Brownian motion with an arbitrary initial law $\alpha$ and $\psi$ is convex, with $\E[|\nabla\psi(B_1)|]<\infty$. 

In dimension one, every Bass martingale is stretched Brownian motion without any additional assumptions~\cite[Theorem~3.1]{BaBeHuKa20}.
In higher dimensions, this equivalence fails.

\begin{theorem}[Main result]\label{thm:main}
    For every $d\ge2$, there exist compactly supported probability measures $\mu\le_{\mathrm{cx}}\nu$ on $\R^d$ with full affine span such that $(\mu,\nu)$ is irreducible and admits a Bass martingale which is not stretched Brownian motion.
    In particular, there exist two Bass martingales with marginals $\mu,\nu$ which have different laws.
\end{theorem}

This corrects the assertion in~\cite[Theorem~1.3]{BaBeScTs26} that every Bass martingale is optimal. The optimality argument invoked there would require additional integrability; for instance, $\E[|B_0|^2]<\infty$ suffices.

Optimality and the Bass property are equivalent if $(\mu, \nu)$ is \emph{strictly irreducible}; see \Cref{cor:strict-monotonicity}.
An irreducible pair $(\mu,\nu)$ is called strictly irreducible if $\operatorname{supp}\mu$ is a compact subset of $\operatorname{ri}{\operatorname{conv}}(\operatorname{supp}\nu)$, cf.~\cite{AcMaPa25, BaPaSc25, CaMaSy26}. 
Indeed, strict irreducibility makes the associated dual potential integrable, so the implication from Bass structure to optimality follows directly from~\cite[Proposition~3.12(2)]{HaJoLoObPa26}.

\subsection*{Cyclical monotonicity and Bass martingales}

Let $\gamma=\gamma_d$ denote the standard Gaussian law on $\R^d$ and $\Pi(p,\gamma)$ the set of couplings of $p$ and $\gamma$.
For $p\in\Pcal_2(\R^d)$, write $\MCov(p,\gamma):=\sup_{q\in\Pi(p,\gamma)}\int y\cdot z\,q(\rmd y,\rmd z)$ for their maximal covariance.
For $\pi\in\Mcal(\mu,\nu)$, let $\pi_x$ denote its conditional terminal law at $x$.
The martingale Benamou--Brenier problem~\eqref{eq:mBB} has the equivalent weak transport formulation~\cite[Theorem~2.2]{BaBeHuKa20}
\begin{equation}\label{eq:mBB-weak} \tag{wmBB}
    \mBB(\mu,\nu)
    =
    \sup_{\pi\in\Mcal(\mu,\nu)}J(\pi),
    \qquad 
    J(\pi):=\int\MCov(\pi_x,\gamma)\,\mu(\rmd x).
\end{equation}
Its unique optimizer, denoted $\operatorname{sBM}(\mu,\nu)$, is the joint law of the endpoints of stretched Brownian motion.
As in classical transport, optimality implies a suitable notion of cyclical monotonicity~\cite{BaBePa18}.

\begin{definition}[MC-cyclical monotonicity]\label{def:mc-monotonicity}
    A martingale coupling $\pi\in\Mcal(\mu,\nu)$ is \emph{MC-cyclically monotone} if there is a $\mu$-full set $G$ such that, for all $x_1,\ldots,x_N\in G$ and $q_1,\ldots,q_N\in\Pcal_2(\R^d)$ with barycenters $x_1,\ldots,x_N$ and $\sum_{i=1}^Nq_i=\sum_{i=1}^N\pi_{x_i}$, one has
    \begin{equation}\label{eq:mc-monotonicity}
        \sum_{i=1}^N\MCov(\pi_{x_i},\gamma)
        \ge
        \sum_{i=1}^N\MCov(q_i,\gamma).
    \end{equation}
\end{definition}
Thus no finite redistribution of the conditional laws can improve the objective while preserving their sum and their barycenters.

A \emph{Bass coupling} is the joint law of the initial and terminal values of a Bass martingale.
As the following theorem shows, every Bass coupling has this property, and the converse holds under the additional \emph{spreading condition} $\pi_x\sim\nu$ for $\mu$-almost every $x$.

\begin{theorem}[Cyclical monotonicity and Bass structure]\label{thm:cyclical-bass}
    Every Bass coupling is MC-cyclically monotone. A martingale coupling satisfying the spreading condition is MC-cyclically monotone if and only if it is a Bass coupling.
\end{theorem}
There is also a geometric sufficient condition for Bass structure: an MC-cyclically monotone coupling is a Bass coupling if $\mu$-almost all starting points lie in the same cell of the convex paving of Ghoussoub--Kim--Lim~\cite{GhKiLi19}; see \Cref{prop:single-cell}.

The spreading condition in \Cref{thm:cyclical-bass} cannot be dropped.
Indeed, \Cref{thm:mc-not-bass} shows that MC-cyclical monotonicity need not imply Bass structure, even for irreducible, compactly supported marginals.
Together with \Cref{thm:main}, this shows that the implications from optimality to Bass structure and from Bass structure to MC-cyclical monotonicity are both strict in general.
Note, however, that under strict irreducibility all three notions are equivalent; see \Cref{cor:strict-monotonicity}.

\subsection*{Duality and Jensen gaps}

All our counterexamples rely on the dependence of Jensen gaps on the martingale coupling.
To explain how this enters the martingale Benamou--Brenier problem, we consider the dual of the weak transport formulation~\eqref{eq:mBB-weak}.
For a convex potential $v$, set $v^C=(v^*\star\gamma)^*$, where ${}^*$ denotes convex conjugation and $\star$ denotes convolution.
Fenchel's inequality gives
\begin{align}
    \MCov(\pi_x,\gamma)
    &\le \int v\,\rmd\pi_x-v^C(x).
    \label{eq:conditional-dual-bound}
\end{align}
Integrating and optimizing over $v$ yields the duality formula~\cite[Theorem~1.4]{BaBeScTs26}
\begin{equation}\label{eq:mBB.dual}\tag{dmBB}
    \mBB(\mu,\nu)
    = \inf_{\substack{v\in L^1(\nu)\\v\text{ convex}}}
      \Bigl\{\int v\,\rmd\nu-\int v^C\,\rmd\mu\Bigr\},
\end{equation}

A dual optimizer need not be $\nu$-integrable.
To accommodate such potentials, one integrates the conditional difference instead of subtracting the marginal integrals.
For convex $v$ with $\mu(\operatorname{dom}v)=1$, the relaxed duality theorem reads
\begin{align}
    \label{eq:mBB.dual.relaxed}\tag{\faUmbrellaBeach-dmBB}
    \mBB(\mu,\nu)
    &=
    \inf_{\substack{v\text{ convex}\\\mu(\operatorname{dom}v)=1}}\sup_{\pi\in\Mcal(\mu,\nu)}\Ecal_\pi(v),
    \\
    \notag
    \Ecal_\pi(v)
    &:=
    \int\Bigl[\int v\,\rmd\pi_x-v^C(x)\Bigr]\,\mu(\rmd x).
\end{align}
By Jensen's inequality and $v^C\le v$, the integrand is nonnegative, so $\Ecal_\pi(v)$ is well defined in $[0,\infty]$ for every $\pi\in\Mcal(\mu,\nu)$.
If $(\mu,\nu)$ is irreducible, the infimum in~\eqref{eq:mBB.dual.relaxed} is attained, and every minimizer $v$ induces a Bass representation of sBM with potential $\psi=v^*$~\cite[Theorem~1.4]{BaBeScTs26}.
Conversely, if $\kappa$ is a Bass coupling with potential $\psi$, the conjugate $v=\psi^*$ gives equality in~\eqref{eq:conditional-dual-bound}; see e.g.~\Cref{prop:Bass.structure}.
If $v\in L^1(\nu)$, this proves optimality, since $\Ecal_\pi(v)=\int v\,\rmd\nu-\int v^C\,\rmd\mu$ depends only on the marginals.
For nonintegrable $v$, however, $\Ecal_\pi(v)$ may depend on $\pi$. 
Equality in~\eqref{eq:conditional-dual-bound} then need not imply optimality.

This dependence is entirely captured by the \emph{Jensen gap}
\begin{equation}\label{eq:jensen-gap}
    D_\pi(v) := \int\Bigl[\int v\,\rmd\pi_x-v(x)\Bigr]\,\mu(\rmd x).
\end{equation}
Indeed, separating the coupling-dependent and coupling-independent terms gives
\begin{equation}\label{eq:dual-jensen-gap}
    \Ecal_\pi(v) = D_\pi(v)+\int\bigl[v(x)-v^C(x)\bigr]\,\mu(\rmd x).
\end{equation}
Conditional differences of this kind also underlie generalized formulations of martingale transport duality~\cite{BeNuTo17,De18,DeTo19}.
For $v\in L^1(\nu)$, the Jensen gap equals $\int v\,\rmd\nu-\int v\,\rmd\mu$ and is therefore independent of $\pi$.
In \Cref{sec:construction}, we show that this independence can fail without integrability in dimensions $d\ge2$, even for irreducible, compactly supported marginals.

\begin{lemma}[Coupling-dependent Jensen gaps]\label{lem:jensen-gaps}
    For every $d\ge2$, there exist an irreducible pair $(\mu,\nu)$ of compactly supported probability measures on $\R^d$ with full affine span, martingale couplings $\pi^1,\pi^2\in\Mcal(\mu,\nu)$, and a nonnegative lower semicontinuous convex function $f:\R^d\to[0,\infty]$ such that
    \[
        \int f\,\rmd\mu=\int f\,\rmd\nu=\infty,
        \qquad
        0=D_{\pi^2}(f)<D_{\pi^1}(f)<\infty.
    \]
    Moreover, $\mu$ is countably supported and $f$ is finite at every atom of $\mu$.
\end{lemma}

In dimension one, by contrast, the Jensen gap is independent of $\pi$ even without integrability~\cite[Appendix~A.3, Lemma~A.15]{BeJu16}; see also~\cite{BeNuTo17}.
Indeed, after subtracting a supporting affine function, every finite convex function on $\R$ can be written as a countable sum of nonnegative Lipschitz convex functions.
Each summand has a coupling-independent Jensen gap, and monotone convergence gives the same conclusion for their sum.
Our example therefore also shows that in dimensions $d\ge2$ not every finite convex function is a countable sum of Lipschitz convex functions; see \Cref{cor:lipschitz}.

\subsection*{Instability of martingale transport}

In dimension one, martingale transport admits general approximation and stability results under convergence of the marginals~\cite{BeJoMaPa22,BaPa22,Wi23}.
In particular, stretched Brownian motion is stable under $\Wcal_2$-convergence of the marginals~\cite[Corollary~1.4]{BeJoMaPa23}.
In dimensions $d\ge2$, however, both properties can fail even when the marginals converge in total variation and are supported on a common compact set.

\begin{theorem}[Instability]\label{thm:instability}
For every $d\ge2$, there exist irreducible pairs $(\mu^n,\nu^n)$ and $(\mu,\nu)$ and a Bass coupling $\kappa\in\Mcal(\mu,\nu)$ such that all marginals have full affine span and are supported in a common compact set, $\mu^n\to\mu$ and $\nu^n\to\nu$ in total variation, and
\begin{align}
    \operatorname{sBM}(\mu^n,\nu^n)
    &\longrightarrow\kappa\ne\operatorname{sBM}(\mu,\nu)
      \quad\text{in total variation},
    \label{eq:instability-couplings}\\
    \mBB(\mu^n,\nu^n)
    &\longrightarrow J(\kappa)<\mBB(\mu,\nu).
    \label{eq:instability-values}
\end{align}
Moreover, no sequence $\pi^n\in\Mcal(\mu^n,\nu^n)$ converges weakly to $\operatorname{sBM}(\mu,\nu)$.
\end{theorem}

This result follows from the existence of a non-optimal Bass coupling established in \Cref{thm:main}.
Keeping its potential fixed, we condition the initial Brownian law on increasing balls. The resulting Bass couplings are stretched Brownian motions, yet they converge in total variation to the original, non-optimal coupling.

Higher-dimensional instability of martingale transport was previously shown by Br\"uckerhoff and Juillet~\cite{BrJu22}.
Our result shows that it persists for the martingale Benamou--Brenier problem even under total-variation convergence, with irreducible marginals supported in a common compact set.

\subsection*{Related literature}

In quadratic optimal transport, the characterization of optimal couplings by convex potentials was established by Brenier~\cite{Br91} and independently by Rachev and R\"uschendorf~\cite{RaRu90}.
McCann's displacement interpolation~\cite{Mc97} and the dynamic formulation of Benamou and Brenier~\cite{BeBr00} connect this geometric structure with a variational principle.
The martingale analogue seeks an interpolation which preserves conditional means rather than following straight lines.

The classical link between cyclical monotonicity and optimality extends far beyond quadratic cost.
Pratelli~\cite{Pr08} treated continuous, possibly infinite-valued costs, and Schachermayer and Teichmann~\cite{ScTe09} treated finite lower semicontinuous costs.
For finite-cost plans, the equivalence holds even for finite Borel costs~\cite[Theorem~1]{BeGoMaSc09}.
Counterexamples require infinite-valued costs: Ambrosio and Pratelli~\cite[Example~3.5]{AmPr03} construct a cost which is finite only on two prescribed graphs.
Our examples instead use the canonical Gaussian weak martingale transport objective.

Martingale optimal transport (MOT) arises naturally in model-independent finance.
For a single asset, European call prices at all strikes determine its marginal law under the pricing measure at each quoted maturity~\cite{BrLi78}, but leave the dynamics between maturities undetermined.
Robust pricing and hedging optimize over compatible martingale models.
This approach was developed through Skorokhod embeddings by Hobson~\cite{Ho98}, and through martingale transport and stochastic control in~\cite{BeHePe13,GaHeTo14,DoSo14,TaTo13}.
Selecting a tractable calibrated interpolation, rather than an extremal price, is a complementary problem for which Bass martingales provide a natural construction.

Bass's original construction solves a Skorokhod embedding problem~\cite{Ba83}.
Ankirchner, Hobson and Strack~\cite{AnHoSt15} extended the martingale-representation approach to embeddings in diffusions.
Related measure-valued constructions were developed by Eldan~\cite{El16}, while Cox and K\"allblad~\cite{CoKa17} used conditional terminal laws as state variables in robust option pricing.

The martingale Benamou--Brenier problem was developed probabilistically, and sBM introduced, in~\cite{BaBeHuKa20}.
Huesmann and Trevisan~\cite{HuTr19} developed a parallel PDE formulation and dynamic duality theory.
The calibration approach of Guo, Loeper and Wang~\cite{GuLoWa19} and Loeper's market-impact model~\cite{Lo18} provide further early formulations of the problem.
The static convex duality and global Bass representation for irreducible marginals were established in~\cite{BaBeScTs26}.
The decomposition of sBM into Bass martingales was subsequently studied by Schachermayer and Tschiderer~\cite{ScTs26}.

A further motivation comes from mathematical finance, where martingale optimal transport arises naturally in model-independent pricing and model calibration.
Dupire's local-volatility construction~\cite{Du94} uses option prices across strikes and maturities to specify a calibrated diffusion.
Variational calibration of local-stochastic and path-dependent volatility models is developed in~\cite{GuLoWa22,GuLo21}.
Conze and Henry-Labord\`ere~\cite{CoHe21} proposed the fixed-point iteration for Bass local-volatility calibration; its existence and convergence theory was developed in~\cite{AcMaPa25}.
The relation between increasingly fine Bass interpolations and the Dupire model is studied in~\cite{BePaSc22,AcMaPa25}.
Variational approaches include the Bass functional~\cite{BaScTs25} and its gradient flow~\cite{BaPaSc25}.
Loeper~\cite{Lo23} formulated the iteration as a martingale Sinkhorn algorithm; its multidimensional formulation and convergence analysis were developed further in~\cite{HaJoLoObPa26}.

The Bass construction has also been extended to positive-price models and other reference laws.
Geometric Bass martingales, introduced by Backhoff-Veraguas, Loeper and Ob{\l}\'oj~\cite{BaLoOb25}, connect the arithmetic construction to positive-price models by a change of num\'eraire~\cite{BePaRi26}.
Extensions beyond Brownian reference dynamics already appear in~\cite[Section~6]{CoHe21}.
Tschiderer's $q$-Bass martingales~\cite{Ts24} allow more general reference laws, with one-dimensional existence and fixed-point results obtained by Acciaio and Marini~\cite{AcMa26a,AcMa26b}.
State-dependent reference laws are considered by Marini~\cite{Ma26qx}.
The stochastic Knothe--Rosenblatt construction~\cite{BeHaPa26} extends this approach to stochastic local volatility (SLV) models, retaining prescribed stochastic factor dynamics while matching option-implied marginals at successive maturities.

Entropy-based calibration leads to martingale Schr\"odinger bridges~\cite{NuWiZh23}.
After suitable rescaling, the Schr\"odinger--Bass bridge interpolates between the classical Schr\"odinger and martingale Benamou--Brenier problems~\cite{HeLoMaPhTo26,HaPaTh26}.

For the proofs of our main results, we use convex duality and the general framework of weak optimal transport introduced by Gozlan, Roberto, Samson and Tetali~\cite{GoRoSaTe17}.
Convexity of dual potentials in weak transport is studied in~\cite{Pr25}, and the boundary behavior of dual optimizing sequences for sBM in~\cite{ScSi25}.

\section{Bass structure, optimality, and stability}\label{sec:structure}

We collect several properties of Bass couplings and their relation to optimality and MC-cyclical monotonicity.
In particular, we show that under the spreading condition, MC-cyclical monotonicity characterizes Bass couplings, while under strict irreducibility all three notions coincide.
We also establish conditions for stability of stretched Brownian motion under perturbations of the marginals.

We begin by recalling the conditional-law representation of Bass couplings and some of its consequences.
Following \cite[Section~3]{BaBeScTs26}, we assume for notational convenience that $\operatorname{aff}(\operatorname{supp}\nu)=\R^d$ whenever considering fixed marginals $(\mu,\nu)$.
Indeed, otherwise we identify this affine hull isometrically with $\R^k$, where $k$ is its dimension, and replace $\gamma$ by a standard Gaussian supported on the affine hull.
This changes $J$ only by an additive constant.

\begin{proposition}[Structure of Bass couplings]\label{prop:Bass.structure}
    Let $\mu,\nu\in\Pcal_2(\R^d)$ and let $\kappa\in\Mcal(\mu,\nu)$ be a Bass coupling with potential $w^*$ and Bass measure $\alpha$.
    Then $\nabla(w^*\star\gamma)$ is a homeomorphism from $\R^d$ onto $\operatorname{int}\operatorname{co}(\operatorname{supp}\nu)$, with inverse $\nabla w^C$, and
    \begin{equation}\label{eq:Bass.structure}
        \kappa_x
        =\Bigl(\nabla w^*\bigl(\cdot+\nabla w^C(x)\bigr)\Bigr)_\#\gamma
        \sim\nu,
        \qquad \mu\text{-a.e. }x.
    \end{equation}
    Moreover,
    \begin{equation}\label{eq:bass-conditional-equality}
        \int w\,\rmd\kappa_x-w^C(x)
        =\MCov(\kappa_x,\gamma),
        \qquad \mu\text{-a.e. }x.
    \end{equation}
    In particular, $\Ecal_\kappa(w)=J(\kappa)$.
\end{proposition}

\begin{proof}
    The Bass representation~\eqref{eq:bass} gives
    \[
        \mu=(\nabla(w^*\star\gamma))_\#\alpha,
        \qquad
        \nu=(\nabla w^*)_\#(\alpha\star\gamma).
    \]
    By the Bass representation, $w^*$ is finite and convex, and $\nu\in\Pcal_2(\R^d)$ implies that its gradient is square-integrable under some translated Gaussian law.
    Moreover, the full affine span of $\nu$ implies $\operatorname{int}(\operatorname{dom}w)\neq\varnothing$.

    According to \cite[Lemma~6.4(ii)--(iii)]{BaBeScTs26}, under these assumptions $w^C=(w^*\star\gamma)^*$ is differentiable on the interior of its domain, and $\nabla(w^*\star\gamma)$ is a homeomorphism onto this interior with inverse $\nabla w^C$.
    Furthermore, \cite[Lemma~6.5]{BaBeScTs26} identifies its range with $\operatorname{int}(\operatorname{dom}w)$.
    Since $\nu=(\nabla w^*)_\#(\alpha\star\gamma)$ and $\alpha\star\gamma$ has a strictly positive density, convex duality gives
    \[
        \operatorname{int}(\operatorname{dom}w)
        =\operatorname{int}\operatorname{co}(\operatorname{supp}\nu).
    \]
    This proves the asserted homeomorphism and, in particular, the differentiability of $w^C$ on $\operatorname{int}\operatorname{co}(\operatorname{supp}\nu)$.
    Moreover, every translated Gaussian law is equivalent to $\alpha\star\gamma$, since both have strictly positive Lebesgue densities.
    Their images under $\nabla w^*$ are therefore equivalent, which gives $\kappa_x\sim\nu$ for $\mu$-a.e.~$x$.

    Fix such $x$ and put $a:=\nabla w^C(x)$.
    The inverse-gradient relation gives $x=\nabla(w^*\star\gamma)(a)$ and $w^C(x)=a\cdot x-\int w^*(a+z)\,\gamma(\rmd z)$.
    Fenchel's equality states that
    \[
        w\bigl(\nabla w^*(a+z)\bigr)+w^*(a+z)
        =
        (a+z)\cdot\nabla w^*(a+z)
        \quad 
        \text{for $\gamma$-a.e.~$z$}.
    \]
    Integrating and using the preceding identities yields
    \[
        \int w\,\rmd\kappa_x-w^C(x)
        =
        \int z\cdot\nabla w^*(a+z)\,\gamma(\rmd z)
        =
        \MCov(\kappa_x,\gamma),
    \]
    where the last equality follows because $z\mapsto\nabla w^*(a+z)$ is the gradient of a convex function and hence an optimal transport map from $\gamma$.
    Integrating against $\mu$ proves $\Ecal_\kappa(w)=J(\kappa)$.
\end{proof}

\paragraph{When Bass martingales are stretched Brownian motion.}
Although every Bass coupling satisfies \eqref{eq:bass-conditional-equality}, it need not maximize $J$, as we show in \Cref{thm:main}.
According to \cite[Proposition~3.12(2)]{HaJoLoObPa26}, for irreducible marginals with finite $p$-moments for some $p>1$, the integrable convex conjugate of a Bass potential attains the dual problem, and the associated Bass coupling is stretched Brownian motion.
In the following, we give several moment conditions that ensure this integrability.

\begin{proposition}\label{lem:integrable-bass}
    Let $\kappa\in\Mcal(\mu,\nu)$ be a Bass coupling with potential $w^*$ and Bass measure $\alpha$.
    If $w\in L^1(\nu)$, then $\kappa=\operatorname{sBM}(\mu,\nu)$.
    In particular, this holds if $\alpha\in\Pcal_p(\R^d)$ and $\nu\in\Pcal_q(\R^d)$ for some $1\le p<\infty$ and $1<q\le\infty$ with $\frac1p+\frac1q=1$.
    This includes the cases where $\alpha\in\Pcal_2(\R^d)$, or where $\alpha\in\Pcal_1(\R^d)$ and $\nu$ is compactly supported.
\end{proposition}

\begin{proof}
    If $w\in L^1(\nu)$, then \eqref{eq:bass-conditional-equality} and the finiteness of $J(\kappa)$ imply $w^C\in L^1(\mu)$.
    Hence, for every $\pi\in\Mcal(\mu,\nu)$, the conditional dual bound \eqref{eq:conditional-dual-bound} gives
    \[
        J(\pi)
        \le \int w\,\rmd\nu-\int w^C\,\rmd\mu
        =J(\kappa).
    \]
    Thus $\kappa$ is stretched Brownian motion.

    It remains to show that the stated moment conditions imply $w\in L^1(\nu)$. Let $B_0\sim\alpha$ and $Z\sim\gamma$ be independent, and set ${M_1:=\nabla w^*(B_0+Z)}$, so that $M_1\sim\nu$.
    Since $w^*$ is finite and convex on $\R^d$, we can choose $y_0\in\partial w^*(0)$.
    Fenchel's inequality at $0$ gives $w(M_1)\ge-w^*(0)$.
    On the other hand, Fenchel's equality at $B_0+Z$ and the subgradient inequality
    $w^*(B_0+Z)\ge w^*(0)+y_0\cdot(B_0+Z)$ give
    \[
        w(M_1)\le (B_0+Z)\cdot(M_1-y_0)-w^*(0).
    \]
    If $\alpha\in\Pcal_p(\R^d)$ and $\nu\in\Pcal_q(\R^d)$ for conjugate exponents $p,q>1$, the upper bound is integrable by H\"older's inequality.
    It is also integrable if $\alpha\in\Pcal_1(\R^d)$ and $\nu$ is compactly supported, since $M_1$ is then bounded.
    The lower bound therefore gives $w\in L^1(\nu)$ in both cases.
\end{proof}

Under strict irreducibility, every Bass coupling is stretched Brownian motion.
Indeed, $\operatorname{supp}\mu$ is a compact subset of $\operatorname{int}\operatorname{conv}(\operatorname{supp}\nu)$, where the inverse-gradient map from \Cref{prop:Bass.structure} is continuous.
Consequently, the Bass measure $\alpha=(\nabla w^C)_\#\mu$ has bounded support, so \Cref{lem:integrable-bass} applies.
Conversely, \cite[Corollary~7.7]{BaBeScTs26} states that every martingale coupling maximizing $J$ between irreducible marginals with finite second moments admits a Bass representation.
Thus Bass structure and optimality coincide under strict irreducibility.

\paragraph{When MC-cyclically monotone couplings have Bass structure.}
We now prove \Cref{thm:cyclical-bass}: every Bass coupling is MC-cyclically monotone, and the converse holds under the spreading condition.

\begin{proof}[Proof of \Cref{thm:cyclical-bass}]
    Let $\kappa$ be a Bass coupling with potential $w^*$, and let $G$ be the $\mu$-full set on which \eqref{eq:bass-conditional-equality} holds.
    For $x_1,\ldots,x_N\in G$ and $q_1,\ldots,q_N$ as in \Cref{def:mc-monotonicity}, Fenchel's inequality and \eqref{eq:bass-conditional-equality} give
    \[
        \sum_{i=1}^N\MCov(q_i,\gamma)
        \le
        \sum_{i=1}^N\Bigl(\int w\,\rmd q_i-w^C(x_i)\Bigr)
        =
        \sum_{i=1}^N\Bigl(\int w\,\rmd\kappa_{x_i}-w^C(x_i)\Bigr)
        =
        \sum_{i=1}^N\MCov(\kappa_{x_i},\gamma).
    \]
    Thus $\kappa$ is MC-cyclically monotone.

    Conversely, suppose $\pi$ is MC-cyclically monotone and $\pi_x\sim\nu$ on a $\mu$-full set $G$, and fix $x_0\in G$.
    For every $x\in G\setminus\{x_0\}$, consider the coupling
    \[
        \chi:=\frac12(\delta_{x_0}\otimes\pi_{x_0}+\delta_x\otimes\pi_x).
    \]
    By MC-cyclical monotonicity, $\chi$ maximizes $J$ among all martingale couplings with the same marginals.
    These marginals are irreducible because $\pi_{x_0}\sim\pi_x$.
    According to \cite[Corollary~7.7]{BaBeScTs26}, every martingale coupling maximizing $J$ for irreducible marginals with finite second moments admits a Bass representation.
    
    Let $w$ be a convex function such that $\nabla w^*$ is the unique convex gradient map transporting $\gamma$ to $\pi_{x_0}$.
    For each $x\in G\setminus\{x_0\}$, the Bass representation of $\chi$ expresses $\pi_{x_0}$ and $\pi_x$ as images of translated Gaussian measures under a common convex gradient map.
    Translating the argument of this map so that it transports $\gamma$ to $\pi_{x_0}$, uniqueness of optimal transport identifies it with $\nabla w^*$.
    Thus $\pi_x$ is the image of a translated Gaussian under $\nabla w^*$, and its barycenter $x$ determines the translation as $\nabla w^C(x)$ by \Cref{prop:Bass.structure}.
    Consequently,
    \[
        \pi_x=\Bigl(\nabla w^*(\cdot+\nabla w^C(x))\Bigr)_\#\gamma,
        \qquad \mu\text{-a.e.~}x.
    \]
    Hence $\pi$ is a Bass coupling with potential $w^*$ and Bass measure $\alpha=(\nabla w^C)_\#\mu$.
\end{proof}

We next give a geometric interpretation of the spreading condition in \Cref{thm:cyclical-bass} in terms of the convex paving introduced by Ghoussoub, Kim and Lim~\cite{GhKiLi19}.

\begin{definition}[Ghoussoub--Kim--Lim convex paving]\label{def:gkl-paving} 
    Let $\Gamma\subseteq\R^d\times\R^d$ be a Borel set, write $\Gamma_x=\{y:(x,y)\in\Gamma\}$, and suppose that $x\in\operatorname{ri}\operatorname{conv}(\Gamma_x)$ whenever $\Gamma_x\neq\varnothing$.
    A \emph{convex paving} of $\Gamma$ is a family $\mathcal C$ of pairwise disjoint relatively open convex sets, each containing at least one point $x$ with $\Gamma_x\neq\varnothing$, such that every such $x$ belongs to a unique cell $C(x)\in\mathcal C$ satisfying
    \[
        \operatorname{ri}\operatorname{conv}(\Gamma_x)\subseteq C(x).
    \]
    The paving is \emph{irreducible} if, for every other convex paving $\mathcal D$ of $\Gamma$, one has $C(x)\subseteq D(x)$ whenever $\Gamma_x\neq\varnothing$.
\end{definition}

According to \cite[Theorem~2.8]{GhKiLi19}, every such $\Gamma$ admits a unique irreducible convex paving.
For MC-cyclically monotone couplings, we show that a cell of full $\mu$-measure implies the spreading condition in \Cref{thm:cyclical-bass}, and hence Bass structure.

\begin{proposition}\label{prop:single-cell}
    Let $\pi\in\Mcal(\mu,\nu)$ be MC-cyclically monotone, and let $G$ be a $\mu$-full Borel set witnessing this property.
    If the irreducible convex paving of
    \[
        \Gamma=\{(x,y):x\in G,\ y\in\operatorname{supp}\pi_x\}
    \]
    has a cell of full $\mu$-measure, then $\pi$ is a Bass coupling.
\end{proposition}

To prove this proposition, we need the following lemma, which allows us to merge conditional laws whose convex supports overlap.
We will also use it in \Cref{sec:construction} to establish irreducibility of the marginals in our Jensen-gap example.

\begin{lemma}\label{lem:merging.conditional.laws}
    Let $\pi\in\Mcal(\mu,\nu)$ and let $a\neq b$ be atoms of $\mu$, with $A:=\mu(\{a\})>0$ and $B:=\mu(\{b\})>0$.
    Suppose that
    \[
        \operatorname{ri}\operatorname{conv}(\operatorname{supp}\pi_a)
        \cap
        \operatorname{ri}\operatorname{conv}(\operatorname{supp}\pi_b)
        \neq\varnothing.
    \]
    Then there exists $\tilde\pi\in\Mcal(\mu,\nu)$ such that
    \[
        \tilde\pi_x=\pi_x
        \quad\mu\text{-a.e. on }\mathbb R^d\setminus\{a,b\},
        \qquad
        \tilde\pi_a\sim\tilde\pi_b\sim\pi_a+\pi_b.
    \]
\end{lemma}

\begin{proof}
    Fix $z$ in the intersection.
    By \cite[Theorem~1.49]{FoSc11}, for every $\rho\in\Pcal_1(\mathbb R^d)$,
    \[
        \operatorname{ri}\operatorname{conv}(\operatorname{supp}\rho)
        =\Bigl\{\int y\,\eta(\rmd y):
        \eta\sim\rho,\;
        \frac{d\eta}{d\rho}\in L^\infty(\rho)\Bigr\}.
    \]
    Thus there exist $r\sim\pi_a$ and $s\sim\pi_b$ with bounded respective densities and common barycenter $z$.
    Choose $h>0$ small enough that
    $\tfrac{h}{A}r\le\tfrac12\pi_a$ and
    $\tfrac{h}{B}s\le\tfrac12\pi_b$, and define
    \[
        p:=\pi_a+\tfrac{h}{A}(s-r),
        \qquad
        q:=\pi_b+\tfrac{h}{B}(r-s).
    \]
    Then
    \[
        p\ge\tfrac12\pi_a+\tfrac{h}{A}s,
        \qquad
        q\ge\tfrac12\pi_b+\tfrac{h}{B}r,
    \]
    so that $p\sim q\sim\pi_a+\pi_b$.
    Moreover,
    \[
        \int y\,p(\rmd y)=a,\qquad
        \int y\,q(\rmd y)=b,\qquad
        Ap+Bq=A\pi_a+B\pi_b.
    \]
    Thus replacing $\pi_a,\pi_b$ by $p,q$ yields the desired martingale coupling.
\end{proof}

\begin{proof}[Proof of \Cref{prop:single-cell}]
    For $x\in G$, let $C_x:=\operatorname{ri}\operatorname{conv}(\operatorname{supp}\pi_x)$.
    Since $x$ is the barycenter of $\pi_x$, we have $x\in C_x$.

    Suppose that $C_x\cap C_{x'}\neq\varnothing$ for distinct $x,x'\in G$, and consider the two-point coupling
    \[
        \chi := \frac12(\delta_x\otimes\pi_x+\delta_{x'}\otimes\pi_{x'}).
    \]
    The construction in the proof of \Cref{lem:merging.conditional.laws} gives another martingale coupling with the same marginals as $\chi$ whose conditional laws are both equivalent to $\frac12(\pi_x+\pi_{x'})$.
    Thus the marginals of $\chi$ are irreducible.
    MC-cyclical monotonicity implies that $\chi$ maximizes $J$ among all martingale couplings between its marginals.
    As in the proof of \Cref{thm:cyclical-bass}, it therefore admits a Bass representation, so \Cref{prop:Bass.structure} gives $\pi_x\sim\pi_{x'}$ and hence $C_x=C_{x'}$.
    
    Thus the distinct sets $C_x$, $x\in G$, form a convex paving of $\Gamma$.
    This paving is irreducible, since every convex paving of $\Gamma$ must contain $C_x$ in the cell containing $x$.  
    If a cell $C$ satisfies $\mu(C)=1$, then $C_x=C$ for all $x\in G\cap C$, so the preceding argument gives $\pi_x\sim\pi_{x'}$ for all $x,x'\in G\cap C$.
    Since $\nu=\int\pi_x\,\mu(\rmd x)$, it follows that $\pi_x\sim\nu$ for $\mu$-almost every $x$.
    The claim follows from \Cref{thm:cyclical-bass}.
\end{proof}

\paragraph{Stability under strict irreducibility.}
For compactly supported marginals $\mu\le_{\mathrm{cx}}\nu$, Carlier, Malamut, and Sylvestre~\cite[Theorem~2.9]{CaMaSy26} prove that strict irreducibility is equivalent to the existence of $\varepsilon>0$ such that
\[
    \Wcal_\infty(\tilde\mu,\mu)\le\varepsilon,
    \qquad
    \int x\,\rmd\tilde\mu=\int x\,\rmd\mu
    \quad\Longrightarrow\quad
    \tilde\mu\le_{\mathrm{cx}}\nu.
\]
Moreover, \cite[Proposition~2.8]{CaMaSy26} shows that this stability implies that every sufficiently small mean-zero $r\in L^\infty(\mu;\R^d)$ can be realized as a conditional barycenter error by some coupling in $\Pi(\mu,\nu)$.
In the following, we extend this correction to simultaneous perturbations of both marginals and use it to establish stability of stretched Brownian motion under strict irreducibility.

\begin{proposition}[Stability under strict irreducibility]\label{prop:stability}
    Suppose $(\mu,\nu)$ is strictly irreducible with $\operatorname{aff}(\operatorname{supp}\nu)=\R^d$.
    Let $\mu^n\le_{\mathrm{cx}}\nu^n$ satisfy
    \[
        \Wcal_\infty(\mu^n,\mu)\to0,
        \qquad
        \Wcal_2(\nu^n,\nu)\to0.
    \]
    Then $\mBB(\mu^n,\nu^n)\to\mBB(\mu,\nu)$, and $\operatorname{sBM}(\mu^n,\nu^n) \to \operatorname{sBM}(\mu,\nu)$ weakly.
\end{proposition}

The proof relies on the following uniform correction of conditional barycenters.
It allows us to restore the martingale constraint for approximate couplings without changing their marginals.

\begin{lemma}\label{lem:uniform-qualification}
    Suppose $(\mu,\nu)$ is strictly irreducible and $\operatorname{aff}(\operatorname{supp}\nu)=\R^d$.
    There exists $\eta>0$ such that, whenever $\mu'\le_{\mathrm{cx}}\nu'$ and
    \[
        \Wcal_\infty(\mu',\mu)+\Wcal_1(\nu',\nu)<\eta,
    \]
    every $r\in L^\infty(\mu';\R^d)$ with $\int r\,\rmd\mu'=0$ and $\|r\|_\infty\le\eta$ admits a coupling $\rho\in\Pi(\mu',\nu')$ satisfying
    \[
        \int y\,\rho_x(\rmd y)=x+r(x),
        \qquad\mu'\text{-a.e.~}x.
    \]
\end{lemma}

\begin{proof}
    We minimally adapt the convexity-gap argument in the proof of \cite[Theorem~2.9]{CaMaSy26}.
    Translate the common barycenter to $0$ and choose $\delta>0$ such that
    \[
        K:=\operatorname{co}(\operatorname{supp}\mu)
            +\overline B(0,\delta)
        \subset\operatorname{int}\operatorname{co}(\operatorname{supp}\nu).
    \]
    The key estimate from that proof is
    \begin{equation}\label{eq:uniform-convex-gap}
        \beta:=
        \inf_{\substack{
            g\text{ convex, $1$-Lipschitz}:\\
            g(0)=0,\ \operatorname{diam}\partial g(K)\ge\tfrac12
        }}
        \Bigl(
            \int g\,\rmd\nu-\int g\,\rmd\mu
        \Bigr)>0.
    \end{equation}
    As shown in the proof of \cite[Theorem~2.9]{CaMaSy26}, the class defining $\beta$ is compact for locally uniform convergence, and the infimum is attained.
    The attained gap is strictly positive, since irreducibility would otherwise force the minimizing $g$ to be affine on $\operatorname{co}(\operatorname{supp}\nu)$, contradicting $\operatorname{diam}\partial g(K)\ge\tfrac12$.

    Choose $0<\eta<\tfrac12\min(\delta,\beta)$ and set $\tilde\mu:=(\operatorname{Id}+r)_\#\mu'$.
    For $\mu',\nu',r$ as in the statement, we have
    \[
        \operatorname{supp}\mu'\cup\operatorname{supp}\tilde\mu\subset K,
        \qquad
        \int x\,\rmd\tilde\mu=\int y\,\rmd\nu',
    \]
    and
    \[
        \Wcal_1(\tilde\mu,\mu)+\Wcal_1(\nu',\nu)
        \le 
        \Wcal_\infty(\mu',\mu)+\|r\|_\infty+\Wcal_1(\nu',\nu)
        <
        2\eta<\beta.
    \]
    
    We claim that $\tilde\mu\le_{\mathrm{cx}}\nu'$. To show this, it suffices to test the convex-order against convex $1$-Lipschitz functions.
    Let $f$ be such a function and set $\Delta:=\operatorname{diam}\partial f(K)$.
    If $\Delta=0$, then $f$ is affine on $K$, and hence
    \[
        \int f\,\rmd\tilde\mu=\int f\,\rmd\mu'\le\int f\,\rmd\nu'.
    \]
    Otherwise, choose a supporting affine function $\ell$ of $f$ at $0$ and define
    \[
        g(x):=\inf_{y\in\R^d}\Bigl\{\frac{f(y)-\ell(y)}{2\Delta}+|x-y|\Bigr\}.
    \]
    Since $0\in K$, all subgradients of $(f-\ell)/(2\Delta)$ on $K$ have norm at most $\tfrac12$, so $g$ agrees with this function and has the same subgradients on $K$.
    Moreover, $g$ is convex and $1$-Lipschitz, with $g(0)=0$ and $\operatorname{diam}\partial g(K)\ge\tfrac12$, and therefore belongs to the class defining $\beta$.
    As $g\le(f-\ell)/(2\Delta)$ everywhere, $\tilde\mu$ is supported on $K$, and the two marginals have the same barycenter,
    \[
        \int f\,\rmd\nu'-\int f\,\rmd\tilde\mu
        \ge
        2\Delta\Bigl(\int g\,\rmd\nu'-\int g\,\rmd\tilde\mu\Bigr)\notag
        \ge
        2\Delta\Bigl(\beta-\Wcal_1(\nu',\nu)-\Wcal_1(\tilde\mu,\mu)\Bigr)\ge0.
    \]
    Thus $\tilde\mu\le_{\mathrm{cx}}\nu'$.

    Strassen's theorem \cite{St65} states that two probability measures with finite first moments are in convex order if and only if they admit a martingale coupling.
    Let $\chi\in\Mcal(\tilde\mu,\nu')$ and set
    \[
        \rho(\rmd x,\rmd y)
        :=
        \mu'(\rmd x)\,\chi_{x+r(x)}(\rmd y).
    \]
    Then $\rho\in\Pi(\mu',\nu')$ and its conditional barycenters are $x+r(x)$, as required.
\end{proof}

\begin{proof}[Proof of \Cref{prop:stability}]
    We first establish
    \begin{equation}\label{ineq:proof:prop:stability:lower.bound}
        \liminf_{n\to\infty}\mBB(\mu^n,\nu^n)\ge\mBB(\mu,\nu).
    \end{equation}
    Let $\pi^{\rm sBM}=\operatorname{sBM}(\mu,\nu)$ and fix $\delta>0$.
    Partition $\operatorname{supp}\mu$ into finitely many Borel sets $A_i$ of positive $\mu$-mass and diameter at most $\delta$, and define $q\in\Pi(\mu,\nu)$ by
    \[
        q_x:=\frac{\pi^{\rm sBM}(A_i\times\cdot)}{\mu(A_i)},
        \qquad x\in A_i.
    \]
    Then $|\int y\,q_x(\rmd y)-x| \le \delta$ and $q_x\le L_\delta\nu$ where $L_\delta:=\max_i\frac1{\mu(A_i)}$, while concavity of $\MCov$ yields $J(q)\ge J(\pi^{\rm sBM})=\mBB(\mu,\nu)$.
    
    By $\Wcal_\infty$-convergence, couple $X_n\sim\mu^n$ and $X\sim\mu$ such that $|X_n-X|\le d_n\to0$ almost surely.
    For each $A_i$, choose a coupling of its common conditional law $q_x$ and $\gamma$ attaining maximal covariance.
    Conditionally on $X\in A_i$, sample $(Y,Z)$ from this coupling, independently of $X_n$ given $X$.
    Then $\operatorname{Law}(Y\mid X)=q_X$, $Z$ is independent of $X_n$, and $\E[Y\cdot Z]=J(q)$.
    Using an optimal $\Wcal_2$-coupling of $\nu$ and $\nu^n$, transport $Y$ to $Y_n\sim\nu^n$, independently of $(X,X_n,Z)$ given $Y$, so that
    \[
        \E[|Y_n-Y|^2]=t_n^2,
        \qquad t_n:=\Wcal_2(\nu^n,\nu)\to0.
    \]
    Let $q^n:=\operatorname{Law}(X_n,Y_n)$ and set $e_n(x):=\int y\,q_x^n(\rmd y)-x$.
    Then
    \[
        |\E[Y-X_n\mid X_n]|\le\delta+d_n,
        \qquad
        \operatorname{Law}(Y\mid X_n)\le L_\delta\nu,
    \]
    and hence
    \[
        \E[|Y_n-Y|\mid X_n]
        \le L_\delta\E[|Y_n-Y|]\le L_\delta t_n.
    \]
    Consequently,
    \[
        \|e_n\|_\infty\le\delta+d_n+L_\delta t_n,
        \qquad
        \int e_n\,\rmd\mu^n=0.
    \]
    Moreover, since $Z$ is independent of $X_n$, Cauchy--Schwarz gives
    \[
        J_{\mu^n}(q^n)
        \ge\E[Y_n\cdot Z]
        \ge\E[Y\cdot Z]-\sqrt d\,t_n
        \ge\mBB(\mu,\nu)-\sqrt d\,t_n.
    \]
    Let $\eta>0$ be the constant from \Cref{lem:uniform-qualification} for $(\mu,\nu)$.
    Choose $\delta>0$ sufficiently small that $\varepsilon:=\sqrt\delta<\min(1,\eta/2)$.
    For all sufficiently large $n$, $\Wcal_\infty(\mu^n,\mu)+\Wcal_1(\nu^n,\nu)<\eta$, and
    \[
        \Bigl\|\frac{1-\varepsilon}{\varepsilon}e_n\Bigr\|_\infty
        \le\frac{\delta+d_n+L_\delta t_n}{\varepsilon}
        \le\eta.
    \]
    Thus \Cref{lem:uniform-qualification} gives $\rho^n\in\Pi(\mu^n,\nu^n)$ with
    \[
        \int y\,\rho_x^n(\rmd y)
        =x-\frac{1-\varepsilon}{\varepsilon}e_n(x),
        \qquad\mu^n\text{-a.e.~}x.
    \]
    Hence $(1-\varepsilon)q^n+\varepsilon\rho^n\in\Mcal(\mu^n,\nu^n)$, and concavity and nonnegativity of $\MCov$ imply
    \[
        \mBB(\mu^n,\nu^n)
        \ge(1-\varepsilon)J_{\mu^n}(q^n)
        \ge(1-\varepsilon)\bigl(\mBB(\mu,\nu)-\sqrt d\,t_n\bigr).
    \]
    Letting first $n\to\infty$ and then $\delta\downarrow0$ proves~\eqref{ineq:proof:prop:stability:lower.bound}.
    
    For the reverse inequality, \cite[Proposition~2.8]{BaBePa18} states that $\pi\mapsto\int C(x,\pi_x)\,\pi(\rmd x,\rmd y)$ is lower semicontinuous in $\Wcal_2$ whenever $C:\Rd\times\Pcal_2(\Rd)\to\R\cup\{+\infty\}$ is bounded below, lower semicontinuous, and convex in its second argument.
    Applying this to $C(x,\rho)=\frac12\int|y|^2\,\rho(\rmd y)-\MCov(\rho,\gamma)$ and using $\Wcal_2(\nu^n,\nu)\to0$, we obtain
    \[
        \limsup_{n\to\infty}\mBB(\mu^n,\nu^n)
        \le J(\pi)\le\mBB(\mu,\nu)
    \]
    for any limit $\pi\in\Mcal(\mu,\nu)$ of optimal couplings along a subsequence realizing the $\limsup$.
    Together with~\eqref{ineq:proof:prop:stability:lower.bound}, this proves convergence of the optimal values, while uniqueness of the optimizer yields weak convergence of $\operatorname{sBM}(\mu^n,\nu^n)$ to $\operatorname{sBM}(\mu,\nu)$.
\end{proof}

If $\operatorname{aff}(\operatorname{supp}\nu)$ is a proper affine subspace of $\R^d$, the same stability result holds provided that $\operatorname{supp}\nu^n\subseteq\operatorname{aff}(\operatorname{supp}\nu)$ for all $n$.
Indeed, convex order then also forces $\mu^n$ to be supported on this affine subspace, so the preceding proof applies verbatim.

With this stability result, we can now show that, under strict irreducibility, MC-cyclical monotonicity, representation as a Bass martingale, and optimality for the martingale Benamou--Brenier problem are equivalent.

\begin{corollary}\label{cor:strict-monotonicity}
    If $(\mu,\nu)$ is strictly irreducible, then MC-cyclical monotonicity, Bass structure, and optimality for the martingale Benamou--Brenier problem are equivalent for couplings in $\Mcal(\mu,\nu)$.
\end{corollary}

\begin{proof}
    It remains to show that MC-cyclical monotonicity implies optimality.
    Let $\pi$ be MC-cyclically monotone, and let $G$ be a $\mu$-full set witnessing this property.
    
    We approximate $\pi$ by sampling.
    For each $n$, partition $\operatorname{supp}\mu$ into finitely many Borel sets $A_i$ of positive $\mu$-mass and diameter at most $\frac1n$.
    For each $i$, take independent samples $(X_{i,j})_j$ from $\mu(\cdot\mid A_i)$ and define
    \[
        \pi^{n,N}:=\sum_i\frac{\mu(A_i)}{N}
        \sum_{j=1}^N\delta_{X_{i,j}}\otimes\pi_{X_{i,j}}.
    \]
    For fixed $n$, the law of large numbers gives $\pi^{n,N}\to\pi$ in $\Wcal_2$ almost surely.
    We may therefore choose $N_n$ and a realization of the samples such that $\Wcal_2(\pi^{n,N_n},\pi)\le \frac1n$ and all sampled points belong to $G$.
    Write $\pi^n:=\pi^{n,N_n}\in\Mcal(\mu^n,\nu^n)$.
    Since $\mu^n(A_i)=\mu(A_i)$ and $\operatorname{diam}(A_i)\le \frac1n$, we also have $\Wcal_\infty(\mu^n,\mu)\le \frac1n$.
    
    In dimension one, \cite[Remark~1.6]{BeJoMaPa23} shows that finite equally weighted MC-cyclically monotone couplings are optimal among martingale couplings with the same marginals.
    The argument extends to arbitrary weights and dimensions, and so $\pi^n = \mathrm{sBM}(\mu^n,\nu^n)$.
    By \Cref{prop:stability}, $\pi^n\to\operatorname{sBM}(\mu,\nu)$ weakly. Since $\pi^n\to\pi$ in $\Wcal_2$, we conclude that $\pi=\operatorname{sBM}(\mu,\nu)$.
\end{proof}

\section{One construction, two Jensen gaps}\label{sec:construction}

In this section, we construct two martingale couplings with the same marginals but different Jensen gaps for a convex function whose integrals under both marginals are infinite.
Starting from an explicit example in $\R^2$, we then obtain irreducibility and compact support before extending the construction to every dimension $d\ge2$.

\begin{lemma}\label{lem:jensen-gaps-2d}
    There exist $\mu,\nu\in\Pcal_1(\R^2)$ with full affine span, $\pi^1,\pi^2\in\Mcal(\mu,\nu)$, and a nonnegative lower semicontinuous convex function $f:\R^2\to[0,\infty]$ such that
    \[
        D_{\pi^1}(f)>D_{\pi^2}(f)=0,
        \qquad
        \int f\,\rmd\mu=\int f\,\rmd\nu=+\infty.
    \]
    Moreover, $\mu$ is countably supported and $f$ is finite at every atom of $\mu$.
\end{lemma}

\begin{proof}[Proof of \Cref{lem:jensen-gaps-2d}]
    For $j\in\mathbb Z$, define
    \[
      \phi_j(s)=2^j s-4^j,\qquad
      \psi_j(t)=2^j t-2\cdot4^j,
      \qquad
      f(s,t)=\sup_j\{\phi_j(s),\psi_j(t)\}.
    \]
    Then $f$ is finite and convex. For $s,t>0$, set $\phi(s):=\max_j\phi_j(s)$ and $\psi(t):=\max_j\psi_j(t)$, so that $f(s,t)=\max\{\phi(s),\psi(t)\}$. For $j\in\mathbb Z$, set
    \[
        \begin{aligned}
          P_j&=2^{j-1}(3,4), & Q_j&=2^{j-1}(4,6),\\
          R_j&=2^{j-1}(3,6), & S_j&=2^{j-1}(6,6).
        \end{aligned}
    \]
    
    We first describe the geometry of $f$; see also \Cref{fig:geometry}. In the first quadrant, the equality curve $\phi=\psi$ consists of the segments $[P_j,Q_j]$ and $[Q_j,P_{j+1}]$, with slopes $2$ and $1$, as shown in the left panel of \Cref{fig:geometry}.
    To identify the active affine pieces of $f$, we compare consecutive terms:
    \[
        \phi_j>\phi_{j-1}\iff s>3\cdot2^{j-1},
        \qquad
        \psi_j>\psi_{j-1}\iff t>3\cdot2^j.
    \]
    Consequently, for $s,t>0$,
    \[
        \phi=\phi_j\iff 3\cdot2^{j-1}\le s\le3\cdot2^j,
        \qquad
        \psi=\psi_j\iff 3\cdot2^j\le t\le3\cdot2^{j+1}.
    \]
    Intersecting these strips with $\{\phi\ge\psi\}$ and $\{\psi\ge\phi\}$, respectively, gives the regions where $f=\phi_j$ and $f=\psi_j$, shown in the right panel of \Cref{fig:geometry}.
    
    \begin{figure}[htbp]
        \centering
    
        \begin{minipage}[t]{0.49\textwidth}
            \centering
            \begin{tikzpicture}[scale=0.53,>=Latex]
                \path[use as bounding box]
                    (0.6,0.25) rectangle (13.4,17.7);
    
                \localframe{a}
                \draw[gray!60,densely dashed,line width=0.55pt]
                    (3,2)--(3,4) (6,2)--(6,8)
                    (2,6)--(4,6) (2,12)--(8,12);
                \draw[black, line width=1.2pt]
                    (2,3)--(3,4)--(4,6)--(6,8)--(8,12)--(12,16);
    
                \begin{scope}
                    \clip (2,2) rectangle (12,16);
                    \targetcircles
                    \sourcecrosses
                \end{scope}
    
                \node[pointlabel,anchor=north west,
                      xshift=1.5pt,yshift=-1.5pt]
                    at (3,3) {$S_{j-1}$};
                \node[pointlabel,anchor=south east,
                      xshift=-1.5pt,yshift=1.5pt]
                    at (3,4) {$P_j$};
                \node[pointlabel,anchor=south east,
                      xshift=-1.5pt,yshift=1.5pt]
                    at (3,6) {$R_j$};
                \node[pointlabel,anchor=north west,
                      xshift=1.5pt,yshift=-1.5pt]
                    at (4,6) {$Q_j$};
                \node[pointlabel,anchor=north west,
                      xshift=1.5pt,yshift=-1.5pt]
                    at (6,6) {$S_j$};
                \node[pointlabel,anchor=south east,
                      xshift=-1.5pt,yshift=1.5pt]
                    at (6,8) {$P_{j+1}$};
                \node[pointlabel,anchor=south east,
                      xshift=-1.5pt,yshift=1.5pt]
                    at (6,12) {$R_{j+1}$};
                \node[pointlabel,anchor=north west,
                      xshift=1.5pt,yshift=-1.5pt]
                    at (8,12) {$Q_{j+1}$};
                \node[pointlabel,anchor=north west,
                      xshift=1.5pt,yshift=-1.5pt]
                    at (12,12) {$S_{j+1}$};
                \node[pointlabel,anchor=south east,
                      xshift=-1.5pt,yshift=1.5pt]
                    at (12,16) {$P_{j+2}$};
    
                \node[pointlabel,anchor=north west,
                      xshift=1.5pt,yshift=-1.5pt]
                    at ({24/7},{30/7}) {$X_{2j-1}$};
                \node[pointlabel,anchor=south east,
                      xshift=-1.5pt,yshift=1.5pt]
                    at ({30/7},{48/7}) {$X_{2j}$};
                \node[pointlabel,anchor=north west,
                      xshift=1.5pt,yshift=-1.5pt]
                    at ({48/7},{60/7}) {$X_{2j+1}$};
                \node[pointlabel,anchor=south east,
                      xshift=-1.5pt,yshift=1.5pt]
                    at ({60/7},{96/7}) {$X_{2j+2}$};
            \end{tikzpicture}
        \end{minipage}
        \hfill
        \begin{minipage}[t]{0.49\textwidth}
            \centering
            \begin{tikzpicture}[scale=0.53,>=Latex]
                \path[use as bounding box]
                    (0.6,0.25) rectangle (13.4,17.7);
    
                \localregions{philow}
                \localframe{a}
    
                \begin{scope}
                    \clip (2,2) rectangle (12,16);
                    \targetcircles
                    \sourcecrosses
                \end{scope}
    
                \localregionlabels{j}{j+1}
            \end{tikzpicture}
        \end{minipage}
    
        \vspace{-1em}
        \caption{For $a=2^{j-1}$, the left panel shows the equality curve $\phi=\psi$, while the right panel shows the regions where $f=\phi_k$ or $f=\psi_k$.
        Crosses and circles mark the atoms of $\mu$ and $\nu$, respectively.}
        \label{fig:geometry}
    \end{figure}
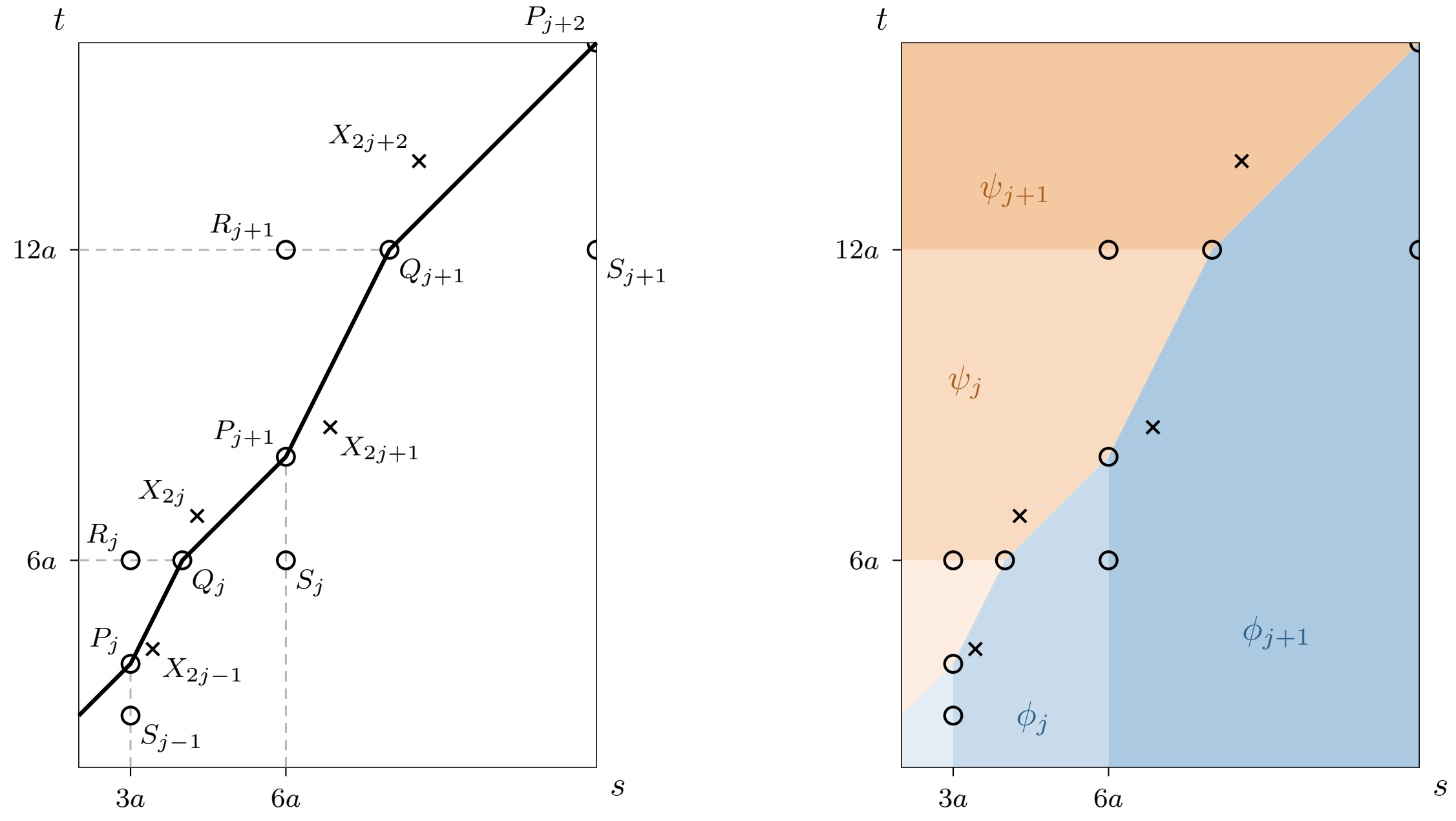

    We next construct the source measure $\mu$, the terminal measure $\nu$, and two martingale couplings $\pi^1,\pi^2 \in \Mcal(\mu,\nu)$. We take the support of $\nu$ to be
    \[
        \{S_0\}\cup\{P_j,Q_j,R_j,S_j:j\ge1\},
        \qquad S_0=(3,3).
    \]
    The support of $\mu$ consists of points $\{X_n\}_{n \ge 0}$ defined as follows: set $X_0=P_1$ and, for $j\ge1$, let $X_{2j-1}$ be the intersection of $[P_j,S_j]$ with $[S_{j-1},Q_j]$, and $X_{2j}$ the intersection of $[Q_j,R_{j+1}]$ with $[R_j,P_{j+1}]$.
    Their coordinates are
    \[
        X_{2j-1}=2^{j-1}\biggl(\frac{24}{7},\frac{30}{7}\biggr),
        \qquad
        X_{2j}=2^{j-1}\biggl(\frac{30}{7},\frac{48}{7}\biggr).
    \]
    By construction, each $X_n$, $n\ge1$, admits two barycentric decompositions in terms of the terminal points:
    \begin{equation}\label{eq:barycenters}
    \begin{aligned}
        X_{2j-1}
        &=\frac67P_j+\frac17S_j
         =\frac47S_{j-1}+\frac37Q_j,\\
        X_{2j}
        &=\frac67Q_j+\frac17R_{j+1}
         =\frac47R_j+\frac37P_{j+1}.
    \end{aligned}
    \end{equation}
    At $X_0=P_1$, we also have $P_1=\frac23S_0+\frac13R_1$.
    We thus define the martingale kernels $(\pi^1_x)_x$ and $(\pi^2_x)_x$ by
    \[
    \begin{aligned}
        \pi^1_{X_0} &= \frac23\delta_{S_0}+\frac13\delta_{R_1},
        & \quad
        \pi^1_{X_{2j-1}} &=\frac67\delta_{P_j}+\frac17\delta_{S_j},
        & \quad
        \pi^1_{X_{2j}} &= \frac67\delta_{Q_j}+\frac17\delta_{R_{j+1}},
        \\
        \pi^2_{X_0} &= \delta_{P_1},
        & \quad
        \pi^2_{X_{2j-1}} &= \frac47\delta_{S_{j-1}}+\frac37\delta_{Q_j},
        & \quad
        \pi^2_{X_{2j}} &= \frac47\delta_{R_j}+\frac37\delta_{P_{j+1}}.
    \end{aligned}
    \]
    The two couplings are illustrated in \Cref{fig:couplings}.
    It remains to choose the masses at the source points $\{X_n\}_{n \ge 0}$ so that both couplings have the same terminal marginal.
    
    \begin{figure}[htbp]
        \centering
    
        \begin{minipage}[t]{0.49\textwidth}
            \centering
            \begin{tikzpicture}[scale=0.53,>=Latex]
                \path[use as bounding box]
                    (0.6,0.25) rectangle (13.4,17.7);
    
                \localregions{excludedcell}
                \localframe{}
                \localregionlabels{1}{2}
    
                \begin{scope}
                    \clip (2,2) rectangle (12,16);
    
                    \draw[firstarrow] (3,4) -- (3,6)
                        node[pos=.57,masssource,
                             text=green!55!black] {$\frac13$};
    
                    \draw[firstarrow] (3,4) -- (3,3)
                        node[pos=.43,masssource,
                             text=green!55!black] {$\frac23$};
    
                    \draw[shortfirstarrow]
                        ({24/7},{30/7})
                        .. controls (4.02,4.18) and (3.92,3.91) ..
                        (3,4)
                        node[pos=.47,below=2.5pt,xshift=-1pt,
                             arrowmass,text=green!55!black]
                        {$\frac67$};
    
                    \draw[firstarrow]
                        ({24/7},{30/7}) -- (6,6);
                    \node[arrowmass, text=green!55!black] at (5.14238,4.77928) {$\frac17$};
    
                    \draw[firstarrow]
                        ({30/7},{48/7}) -- (4,6);
                    \node[arrowmass, text=green!55!black] at (3.71219,6.59117) {$\frac67$};
    
                    \draw[firstarrow]
                        ({30/7},{48/7}) -- (6,12);
                    \node[arrowmass, text=green!55!black] at (4.87239,10.26158) {$\frac17$};
    
                    \draw[firstarrow]
                        ({48/7},{60/7}) -- (6,8);
                    \node[arrowmass, text=green!55!black] at (6.69476,7.88637) {$\frac67$};
    
                    \draw[firstarrow]
                        ({48/7},{60/7}) -- (12,12);
                    \node[arrowmass, text=green!55!black] at (9.97415,10.02448) {$\frac17$};
    
                    \draw[firstarrow]
                        ({60/7},{96/7}) -- (8,12);
                    \node[arrowmass, text=green!55!black] at (7.83035,13.00891) {$\frac67$};
    
                    \draw[firstarrow]
                        ({60/7},{96/7}) -- (12,24);
                    \node[arrowmass, text=green!55!black] at (8.38486,14.86218) {$\frac17$};
    
                    \targetcircles
                    \sourcecrosses
                    \crossat{(3,4)}
                \end{scope}
            \end{tikzpicture}
        \end{minipage}
        \hfill
        \begin{minipage}[t]{0.49\textwidth}
            \centering
            \begin{tikzpicture}[scale=0.53,>=Latex]
                \path[use as bounding box]
                    (0.6,0.25) rectangle (13.4,17.7);
    
                \localregions{excludedcell}
                \localframe{}
                \localregionlabels{1}{2}
    
                \begin{scope}
                    \clip (2,2) rectangle (12,16);
    
                    \draw[secondarrow,shorten <=0pt,shorten >=0pt]
                        (2.83,4.17)
                        .. controls (1.95,4.93) and (1.95,3.10) ..
                        (2.83,3.83);
                    \node[arrowmass,text=violet!85!black]
                        at (2.45,4.7) {$1$};
    
                    \draw[secondarrow]
                        ({24/7},{30/7}) -- (3,3);
                    \node[arrowmass, text=violet!85!black] at (3.7143,3.5619) {$\frac47$};
    
                    \draw[secondarrow]
                        ({24/7},{30/7}) -- (4,6);
                    \node[arrowmass, text=violet!85!black] at (4.21143,5.05331) {$\frac37$};
    
                    \draw[secondarrow]
                        ({30/7},{48/7}) -- (3,6);
                    \node[arrowmass, text=violet!85!black] at (3.38774,6.8113) {$\frac47$};
    
                    \draw[secondarrow]
                        ({30/7},{48/7}) -- (6,8);
                    \node[arrowmass, text=violet!85!black] at (4.9452,7.87357) {$\frac37$};
    
                    \draw[secondarrow]
                        ({48/7},{60/7}) -- (6,6);
                    \node[arrowmass, text=violet!85!black] at (6.88394,7.13394) {$\frac47$};
    
                    \draw[secondarrow]
                        ({48/7},{60/7}) -- (8,12);
                    \node[arrowmass, text=violet!85!black] at (7.92578,10.19617) {$\frac37$};
    
                    \draw[secondarrow]
                        ({60/7},{96/7}) -- (6,12);
                    \node[arrowmass, text=violet!85!black] at (7.01947,13.25652) {$\frac47$};
    
                    \draw[secondarrow]
                        ({60/7},{96/7}) -- (12,16);
                    \node[arrowmass, text=violet!85!black] at (9.80264,15.13603) {$\frac37$};
    
                    \targetcircles
                    \sourcecrosses
                    \crossat{(3,4)}
                \end{scope}
            \end{tikzpicture}
        \end{minipage}
    
        \vspace{-1em}
        \caption{The martingale couplings $\pi^1$ (left, green) and $\pi^2$ (right, violet), shown against the active affine regions of $f$.
        Crosses and circles mark the atoms of $\mu$ and $\nu$, respectively.
        The fractions give the proportions in which mass is transported from each source point to the terminal points.}
        \label{fig:couplings}
    \end{figure}
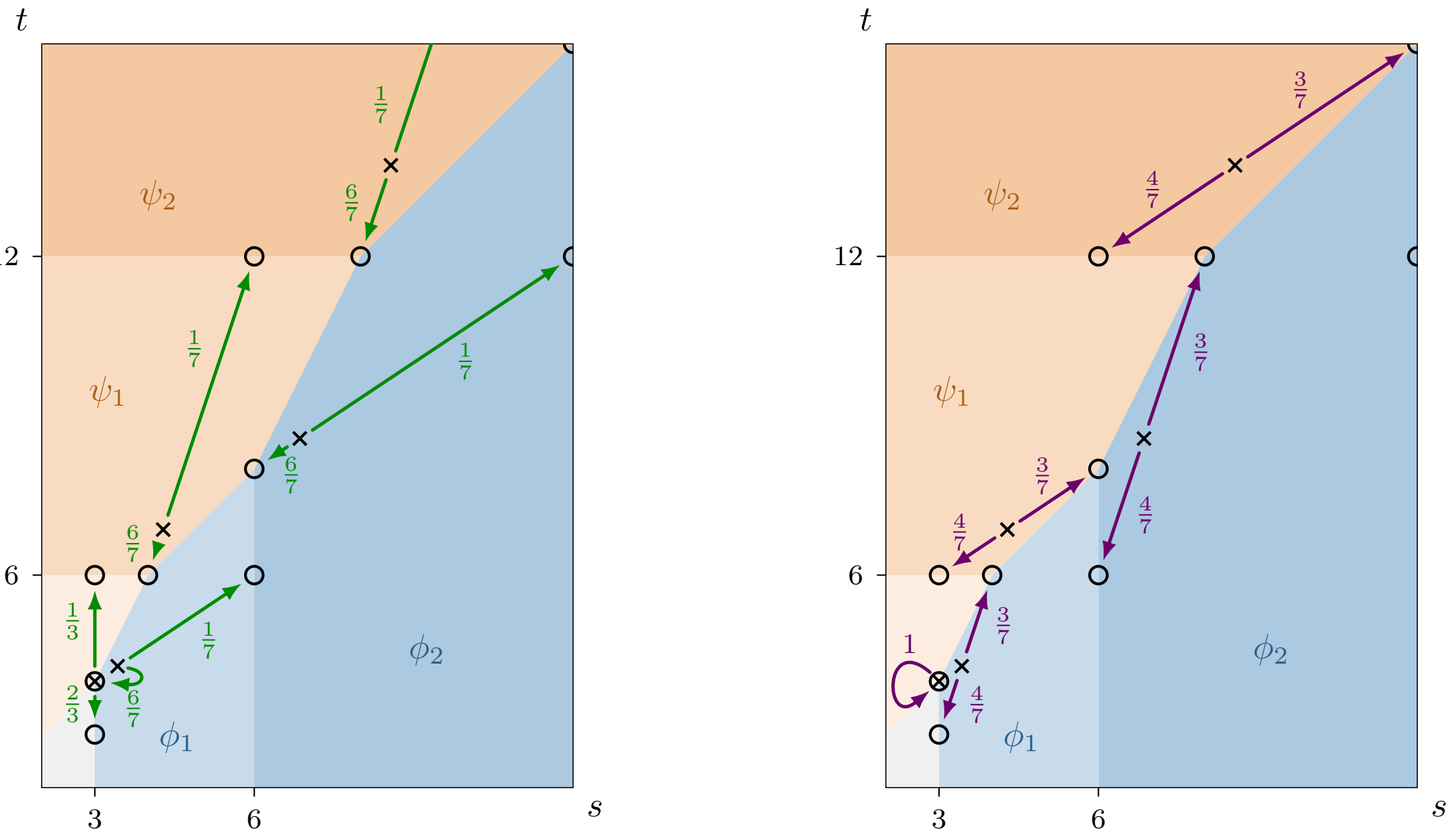

    Write $\mu_n=\mu(\{X_n\})$.
    For the two couplings to have the same terminal marginal, the masses arriving at $Q_j$ and $P_{j+1}$ must agree.
    This gives
    \[
        \frac37\mu_{2j-1}=\frac67\mu_{2j},
        \qquad
        \frac37\mu_{2j}=\frac67\mu_{2j+1}.
    \]
    Hence $\mu_{n+1}=\mu_n/2$ for $n\ge1$, while equality at $P_1$ requires $\mu_0=\frac67\mu_1$.
    Together with $\sum_{n\ge0}\mu_n=1$, these relations yield
    \[
        \mu_0=\frac3{10},
        \qquad
        \mu_n=\frac7{20}\cdot2^{-(n-1)},\quad n\ge1.
    \]
    The masses also agree at $S_0$ and $R_1$, and at all remaining $S_j$ and $R_j$ by the relation $\mu_{n+2}=\frac14\mu_n$.
    
    Thus $\pi^1,\pi^2$ have common marginals $\mu=\sum_{n\ge0}\mu_n\delta_{X_n}$ and $\nu$.
    For $n\ge1$, we have
    \[
        \mu_{n+2}=\frac14\mu_n,
        \qquad
        X_{n+2}=2X_n.
    \]
    Since the terminal points grow on the same dyadic scale, both marginals have finite first moments.
    Finally, the identities \eqref{eq:barycenters}, together with the decomposition at $X_0$, show that $\pi^1$ and $\pi^2$ are martingale couplings of $\mu$ and $\nu$. In particular, $\mu\le_{\mathrm{cx}}\nu$.
    
    For $n\ge1$, both couplings transport mass along segments on which $f$ is affine.
    Indeed, $f$ agrees with $\phi_j$ on both segments through $X_{2j-1}$ and with $\psi_j$ on both segments through $X_{2j}$. Consequently, for every $n\ge1$,
    \[
        \int f\,\rmd\pi^1_{X_n}-f(X_n) = \int f\,\rmd\pi^2_{X_n}-f(X_n) = 0.
    \]
    On the other hand, at $X_0=P_1$, the coupling $\pi^2$ leaves the mass at $P_1$, whereas $\pi^1$ splits it between $S_0$ and $R_1$.
    Since $f(S_0)=f(P_1)=2$ and $f(R_1)=4$, we have
    \[
        \int f\,\rmd\pi^1_{X_0}-f(X_0)
        =
        \frac23,
        \qquad
        \int f\,\rmd\pi^2_{X_0}-f(X_0)
        =
        0.
    \]
    Hence
    \[
        D_{\pi^1}(f)=\frac23\mu_0=\frac15,
        \qquad
        D_{\pi^2}(f)=0.
    \]
    Finally, $\int f\,\rmd\mu=\int f\,\rmd\nu=\infty$, since $\mu_n f(X_n)=1$ for every $n\ge1$, and Jensen's inequality gives $\int f\,\rmd\nu\ge\int f\,\rmd\mu$.
\end{proof}

In dimension one, the Jensen gap depends only on the marginals, even without integrability.
The proof of~\cite[Lemma~A.15]{BeJu16} approximates convex functions by extending them affinely outside increasing bounded intervals.
These approximations are Lipschitz, and their successive differences are again convex, allowing one to pass to the limit in the Jensen gaps by monotone convergence, even without integrability; see also~\cite{BeNuTo17}.
In particular, every finite convex function on $\mathbb R$ can be written as a pointwise convergent sum of Lipschitz convex functions.
Our construction shows that this is no longer true in dimensions $d\ge2$.

\begin{corollary}\label{cor:lipschitz}
    For every $d\ge2$, there exists a finite convex function on $\mathbb R^d$ which cannot be written as a pointwise convergent countable sum of Lipschitz convex functions.
\end{corollary}

\begin{proof}
    Let $f$ be the function from \Cref{lem:jensen-gaps-2d}, and suppose that $f=\sum_{k\ge1}g_k$ pointwise, where each $g_k$ is Lipschitz and convex.
    Since $0\le f(z)\le |z|^2/4$, we have $f(0)=0$ and $\nabla f(0)=0$.
    For each $k$, choose a supporting affine function $\ell_k$ of $g_k$ at the origin and set $h_k=g_k-\ell_k\ge0$.
    The convex partial sums $\sum_{k=1}^N g_k$ converge locally uniformly to $f$.
    Since $f$ is differentiable at the origin, their supporting affine functions $\sum_{k=1}^N\ell_k$ converge to zero, and hence $f=\sum_{k\ge1}h_k$.
    
    For every martingale coupling $\pi \in \Mcal(\mu,\nu)$, Tonelli's theorem gives
    \[
        \int\Bigl[\int f\,\rmd\pi_x - f(x)\Bigr]\,\mu(\rmd x)
        =
        \sum_k\int\Bigl[\int h_k \,\rmd \pi_x-h_k(x)\Bigr]\,\mu(\rmd x)
        =
        \sum_k \Bigl(\int h_k\,\rmd\nu-\int h_k\,\rmd\mu\Bigr).
    \]
    Each term is finite because $h_k$ is Lipschitz and the marginals
    have finite first moments. The right-hand side is independent of $\pi$, contradicting \Cref{lem:jensen-gaps-2d}.
\end{proof}

Using \Cref{lem:merging.conditional.laws}, we modify the example from \Cref{lem:jensen-gaps-2d} to obtain irreducible marginals while preserving the strict inequality between the Jensen gaps.

\begin{lemma}\label{lem:irreducible.jensen-gaps}
    There exist an irreducible pair $(\mu,\nu)$ of probability measures in $\Pcal_1(\R^2)$ with full affine span, $\pi^1,\pi^2\in\Mcal(\mu,\nu)$, and a nonnegative lower semicontinuous convex function $f:\R^2\to[0,\infty]$ such that
    \[
        D_{\pi^1}(f)>D_{\pi^2}(f)=0,
        \qquad
        \int f\,\rmd\mu=\int f\,\rmd\nu=+\infty.
    \]
    Moreover, $\mu$ is countably supported and $f$ is finite at every atom of $\mu$.
\end{lemma}

\begin{proof}
    Let $\mu,\nu,f,\pi^1,\pi^2$ be as in \Cref{lem:jensen-gaps-2d}.
    Set $Q_I=(0,\infty)^2$ and let, for $i \ge 1$ and $j \ge 0$,
    \[
        \Phi_i:=\operatorname{int}\{(s,t)\in Q_I:f(s,t)=\phi_i(s)\},
        \qquad
        \Psi_j:=\operatorname{int}\{(s,t)\in Q_I:f(s,t)=\psi_j(t)\}.
    \]
    These cells are nonempty, bounded and convex, with $\Phi_i=2^{i-1}\Phi_1$ and $\Psi_j=2^j\Psi_0$.
    Therefore $A_i:=2^{i-1}(5,1)\in \Phi_i$ and $B_j:=2^j(1,5)\in \Psi_j$.
    These points are distinct from the atoms $\{X_n\}_n$ of $\mu$ and $\{P_n,Q_n,R_n,S_n\}_n$ of $\nu$.
    By \cite[Theorem~1.49]{FoSc11}, every $\rho\in\Pcal_1(\R^2)$ satisfies
    \[
        \operatorname{ri}\operatorname{conv}(\operatorname{supp}\rho)
        =
        \Bigl\{\int y\,\eta(\rmd y): \eta\sim\rho,\; \frac{\rmd\eta}{\rmd\rho}\in L^\infty(\rho)\Bigr\}.
    \]
    Since $C=\operatorname{ri}\operatorname{conv}(\operatorname{supp}U_C)$ for the uniform law $U_C$ on each cell $C$, we may choose $\theta_{\Phi_i}\sim U_{\Phi_i}$ and $\theta_{\Psi_j}\sim U_{\Psi_j}$ with barycenters $A_i$ and $B_j$, respectively.

    Set $\mathcal C:=\{\Phi_i\}_i\cup\{\Psi_j\}_j$ and choose $w_C>0$ with $\sum_Cw_C=1$ and $\sum_C w_C \sup_{y\in C}|y|<\infty$.
    Define
    \[
        \xi:=\sum_Cw_C\,\delta_{\bar\theta_C}\otimes\theta_C,
        \qquad
        \tilde\pi^k:=\tfrac12(\pi^k+\xi),\quad k=1,2.
    \]
    Then $\tilde\pi^1,\tilde\pi^2\in\Mcal(\tilde\mu,\tilde\nu)$, where
    \[
        \tilde\mu=\tfrac12\mu+\tfrac12\sum_Cw_C\delta_{\bar\theta_C},
        \qquad
        \tilde\nu=\tfrac12\nu+\tfrac12\sum_Cw_C\theta_C.
    \]
    Since $f$ is affine on each cell, $D_\xi(f)=0$, and hence $D_{\tilde\pi^1}(f)=\tfrac12D_{\pi^1}(f)>0=D_{\tilde\pi^2}(f)$.
    
    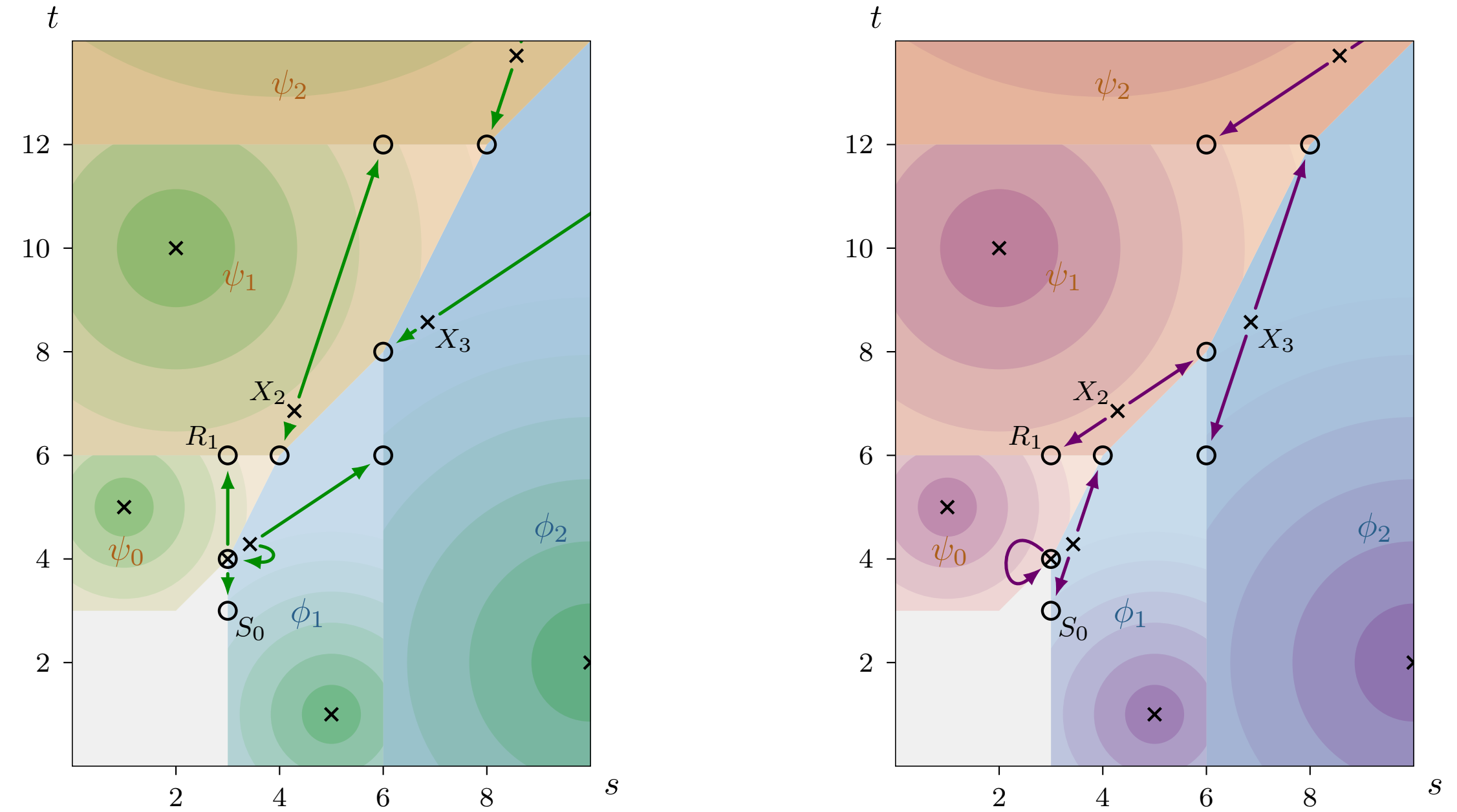
\begin{figure}[t]
        \centering
    
        \begin{minipage}[t]{0.49\textwidth}
            \centering
            \begin{tikzpicture}[scale=0.53,>=Latex]
                \path[use as bounding box]
                    (-0.8,-0.75) rectangle (12.0,16.7);
    
                \constructionregions
    
                \begin{scope}
                    \clip (0,0) rectangle (10,14);
                    \addedhalos{green!55!black}
    
                    \draw[firstarrow] (3,4) -- (3,6);
                    \draw[firstarrow] (3,4) -- (3,3);
    
                    \draw[shortfirstarrow]
                        ({24/7},{30/7})
                        .. controls (4.02,4.18) and (3.92,3.91) ..
                        (3,4);
    
                    \draw[firstarrow]
                        ({24/7},{30/7}) -- (6,6);
                    \draw[firstarrow]
                        ({30/7},{48/7}) -- (4,6);
                    \draw[firstarrow]
                        ({30/7},{48/7}) -- (6,12);
                    \draw[firstarrow]
                        ({48/7},{60/7}) -- (6,8);
                    \draw[firstarrow]
                        ({48/7},{60/7}) -- (12,12);
                    \draw[firstarrow]
                        ({60/7},{96/7}) -- (8,12);
                    \draw[firstarrow]
                        ({60/7},{96/7}) -- (12,24);
                \end{scope}
    
                \constructionframe
                \constructionlabels
    
                \begin{scope}
                    \clip (0,0) rectangle (10,14);
                    \originalmarks
                    \atomlabels
                \end{scope}
            \end{tikzpicture}
        \end{minipage}
        \hfill
        \begin{minipage}[t]{0.49\textwidth}
            \centering
            \begin{tikzpicture}[scale=0.53,>=Latex]
                \path[use as bounding box]
                    (-0.8,-0.75) rectangle (12.0,16.7);
    
                \constructionregions
    
                \begin{scope}
                    \clip (0,0) rectangle (10,14);
                    \addedhalos{violet!85!black}
    
                    \draw[secondarrow,shorten <=0pt,shorten >=0pt]
                        (2.83,4.17)
                        .. controls (1.95,4.93) and (1.95,3.10) ..
                        (2.83,3.83);
    
                    \draw[secondarrow]
                        ({24/7},{30/7}) -- (3,3);
                    \draw[secondarrow]
                        ({24/7},{30/7}) -- (4,6);
                    \draw[secondarrow]
                        ({30/7},{48/7}) -- (3,6);
                    \draw[secondarrow]
                        ({30/7},{48/7}) -- (6,8);
                    \draw[secondarrow]
                        ({48/7},{60/7}) -- (6,6);
                    \draw[secondarrow]
                        ({48/7},{60/7}) -- (8,12);
                    \draw[secondarrow]
                        ({60/7},{96/7}) -- (6,12);
                    \draw[secondarrow]
                        ({60/7},{96/7}) -- (12,16);
                \end{scope}
    
                \constructionframe
                \constructionlabels
    
                \begin{scope}
                    \clip (0,0) rectangle (10,14);
                    \originalmarks
                    \atomlabels
                \end{scope}
            \end{tikzpicture}
        \end{minipage}
    
        \vspace{-.3em}
        \caption{The martingale couplings $\tilde\pi^1$ (left) and $\tilde\pi^2$ (right).
        Crosses denote source atoms and circles denote target atoms.
        The concentric shading illustrates transport from the source atoms $A_i$ and $B_j$ into the corresponding affine regions $\Phi_i$ and $\Psi_j$.}
    \end{figure}

    It remains to show that $(\tilde\mu,\tilde\nu)$ is irreducible.
    Starting from $\tilde\pi^1$, we repeatedly use \Cref{lem:merging.conditional.laws} to \emph{merge} each new conditional law with those previously incorporated, making all equivalent to their sum while preserving the marginals and martingale property.    
    We first merge the conditional law supported on $[S_0,R_1]$ with $\theta_{\Psi_0}$, since $(3,5)\in\operatorname{ri}[S_0,R_1]\cap \Psi_0$.
    We then incorporate $\theta_{\Phi_1}$, since
    \[
        \Phi_1\cap\operatorname{int}\operatorname{conv}([S_0,R_1]\cup \Psi_0) \neq \varnothing.
    \]
    More generally, the geometry of the cells gives, for every $j\ge1$,
    \[
        2^{j-1}X_2\in \Psi_j\cap\operatorname{int}\operatorname{conv}(\Phi_j\cup \Psi_{j-1}),
        \qquad
        2^{j-1}X_3\in \Phi_{j+1}\cap\operatorname{int}\operatorname{conv}(\Phi_j\cup \Psi_j).
    \]
    Indeed, these inclusions follow from those for $j=1$ as well as the identities $\Phi_j=2^{j-1}\Phi_1$ and $\Psi_j=2^j\Psi_0$.
    Thus, having incorporated $\theta_{\Phi_j}$ and $\theta_{\Psi_{j-1}}$, we can next merge with $\theta_{\Psi_j}$ and then with $\theta_{\Phi_{j+1}}$.
    Repeating this incorporates all auxiliary conditional laws.
    Whenever $\theta_{\Phi_j}$ or $\theta_{\Psi_j}$ is incorporated, we also merge with $\tilde\pi^1_{X_{2j-1}}$ or $\tilde\pi^1_{X_{2j}}$, respectively.
    This is possible since $X_{2j-1}\in \Phi_j$ and $X_{2j}\in \Psi_j$, and each $X_k$ belongs to the relative interior of its conditional convex support.
    
    Enumerating the source atoms in this order as $\{z_i\}_{i\ge1}$, we obtain $\chi^n\in\Mcal(\tilde\mu,\tilde\nu)$ for each $n$ satisfying
    \[
        \chi^n_{z_i}\sim\sum_{j=1}^n\tilde\pi^1_{z_j},
        \qquad i\le n.
    \]
    Since $\chi^n$ differs from $\tilde\pi^1$ at only finitely many atoms, where its conditional laws have bounded support, we have $D_{\chi^n}(f)<\infty$.
    Choose $a_n>0$ with $\sum_{n\ge1}a_n=1$ and $\sum_{n\ge1}a_nD_{\chi^n}(f)<\infty$, and set $\chi:=\sum_{n\ge1}a_n\chi^n\in\Mcal(\tilde\mu,\tilde\nu)$.
    Then 
    \[
        D_\chi(f)=\sum_{n\ge1}a_nD_{\chi^n}(f)<\infty.
    \]
    For any $i,j$, choosing $n\ge\max\{i,j\}$ gives $\tilde\pi^1_{z_j}\ll\chi^n_{z_i}\ll\chi_{z_i}$, hence $\tilde\nu\ll\chi_{z_i}$.
    Conversely, $\chi_{z_i}\ll\tilde\nu$ by the common second marginal.
    Thus $\chi_{z_i}\sim\tilde\nu$ for every $i$, proving irreducibility.
\end{proof}

The change of num\'eraire is a classical tool in mathematical finance \cite{GeElRo95}, applied to martingale optimal transport by \cite{CaLaMa17} and extended to weak martingale transport by \cite{BePaRi26}.
We use it to compactify the preceding construction while preserving irreducibility and the Jensen gaps.

\begin{lemma}\label{lem:compact.irreducible.jensen-gaps}
    There exist an irreducible pair $(\mu,\nu)$ of compactly supported probability measures on $\R^2$ with full affine span, $\pi^1,\pi^2\in\Mcal(\mu,\nu)$, and a nonnegative lower semicontinuous convex function $f:\R^2\to[0,\infty]$ such that
    \[
        D_{\pi^1}(f)>D_{\pi^2}(f)=0,
        \qquad
        \int f\,\rmd\mu=\int f\,\rmd\nu=+\infty.
    \]
    Moreover, $\mu$ is countably supported and $f$ is finite at every atom of $\mu$.
\end{lemma}

\begin{proof}
    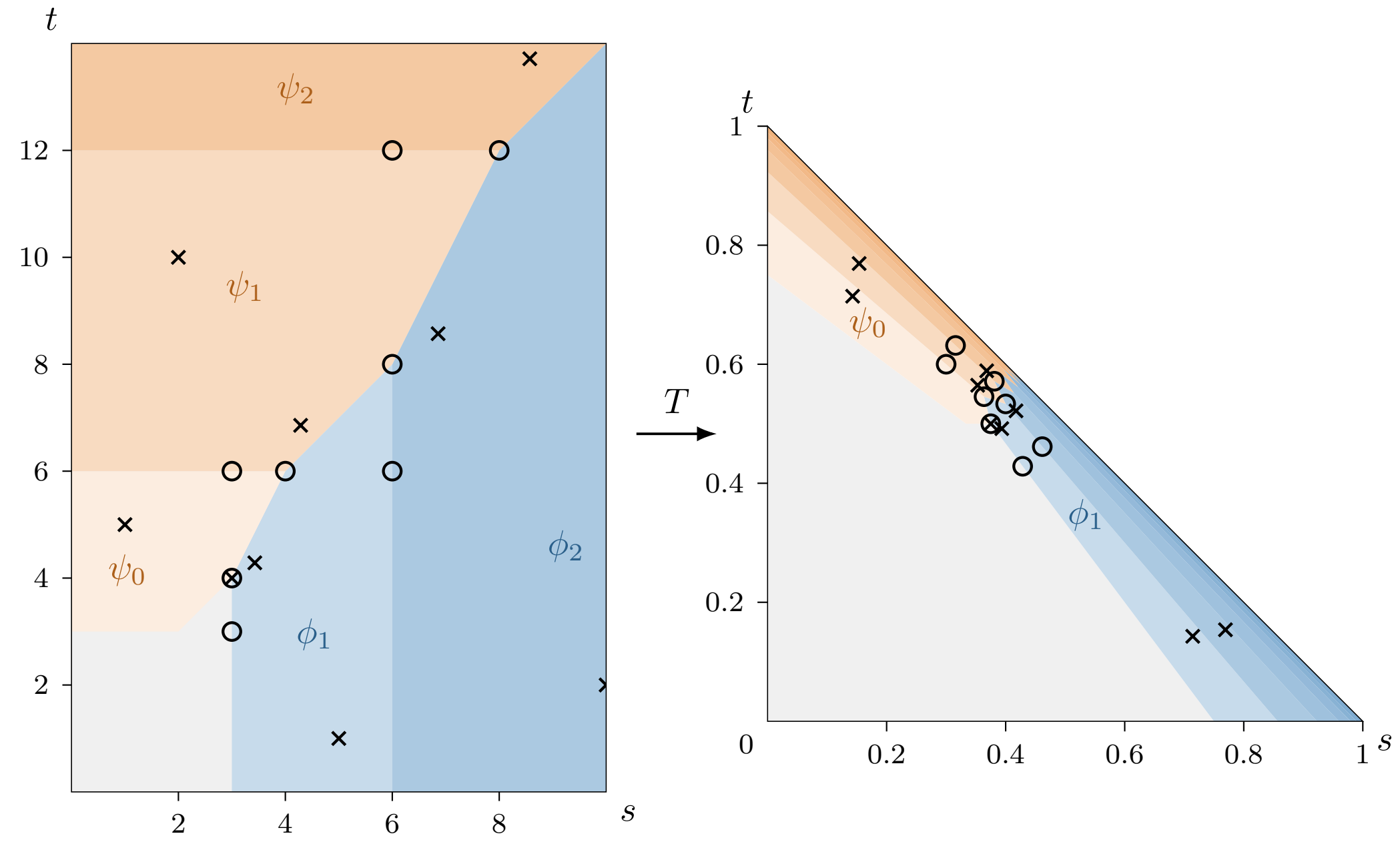
\begin{figure}[htbp]
        \centering
    
        \begin{tikzpicture}[>=Latex]
            \path[use as bounding box]
                (0,0) rectangle (13.9,8.6);

            \begin{scope}[shift={(0.7,0.65)},scale=0.53]
                \constructionregions
    
                \begin{scope}
                    \clip (0,0) rectangle (10,14);
                    \originalmarks
                \end{scope}
    
                \constructionframe
                \constructionlabels
            \end{scope}

            \draw[-{Latex[length=2mm]},black,line width=.72pt]
                (6.30,4.20)--(7.10,4.20);
            \node[above,font=\small,text=black]
                at (6.70,4.27) {$T$};

            \begin{scope}[shift={(7.60,1.35)},x=5.9cm,y=5.9cm]
                \begin{scope}
                    \clip (0,0)--(1,0)--(0,1)--cycle;

                    \fill[excludedcell]
                        \Tpt{0}{0}--\Tpt{3}{0}--
                        \Tpt{3}{4}--\Tpt{2}{3}--
                        \Tpt{0}{3}--cycle;

                    \foreach \i in {1,...,10}{
                        \compactphicell{\i}
                    }
                    \foreach \j in {0,...,9}{
                        \compactpsicell{\j}
                    }
                    \Tcrossat{3}{4}
                    \Tcircleat{3}{4}
                    \Tcrossat{24/7}{30/7}
                    \Tcrossat{30/7}{48/7}
                    \Tcrossat{5}{1}
                    \Tcrossat{1}{5}
                    \Tcircleat{3}{3}
                    \Tcircleat{3}{6}
                    \Tcircleat{4}{6}
                    \Tcircleat{6}{8}
                    \Tcrossat{48/7}{60/7}
                    \Tcrossat{60/7}{96/7}
                    \Tcrossat{10}{2}
                    \Tcrossat{2}{10}
                    \Tcircleat{6}{6}
                    \Tcircleat{6}{12}
                    \Tcircleat{8}{12}
                \end{scope}
    
                \draw[plotframe]
                    (0,0)--(1,0)--(0,1)--cycle;
    
                \foreach \xx/\lab in {
                    0.2/0.2,0.4/0.4,0.6/0.6,
                    0.8/0.8,1/1
                } {
                    \draw[black] (\xx,0)--++(0,-.017)
                        node[below,font=\scriptsize] {$\lab$};
                }
    
                \foreach \yy/\lab in {
                    0.2/0.2,0.4/0.4,0.6/0.6,
                    0.8/0.8,1/1
                } {
                    \draw[black] (0,\yy)--++(-.017,0)
                        node[left,font=\scriptsize] {$\lab$};
                }
    
                \node[below left,font=\scriptsize]
                    at (0,0) {$0$};
                \node[below right,font=\small]
                    at (1,0) {$s$};
                \node[above left,font=\small]
                    at (0,1) {$t$};

                \node[philabel]
                    at \Tpt{4.55}{2.95} {$\phi_1$};
                \node[psilabel]
                    at \Tpt{1.05}{4.15} {$\psi_0$};
            \end{scope}
        \end{tikzpicture}
    
        \caption{Active affine regions of the function $f$ before and after the compactification $T(s,t)=(s,t)/(1+s+t)$.
        Darker shades indicate higher-index cells; crosses and circles mark source and target atoms.}
        \label{fig:compactification-regions}
    \end{figure}

    Let $\mu,\nu,\pi^1,\pi^2,f$ be as in \Cref{lem:irreducible.jensen-gaps}.
    Set $h(s,t):=1+s+t$, $T(s,t):=(s,t)/h(s,t)$ and let $c:=\int h\,\rmd\mu=\int h\,\rmd\nu$.
    Define
    \[
        \hat\mu(\rmd x):=\frac{h(x)}{c}\mu(\rmd x),
        \qquad
        \hat\nu(\rmd y):=\frac{h(y)}{c}\nu(\rmd y),
        \qquad
        \hat\pi^k(\rmd x,\rmd y):=\frac{h(y)}{c}\pi^k(\rmd x,\rmd y),
    \]
    and set
    \[
        \tilde\mu:=T_\#\hat\mu,\qquad
        \tilde\nu:=T_\#\hat\nu,\qquad
        \tilde\pi^k:=(T,T)_\#\hat\pi^k,\quad k=1,2.
    \]
    Since $\int h(y)\,\pi^k_x(\rmd y)=h(x)$, the marginals of $\hat\pi^k$ are $\hat\mu$ and $\hat\nu$.
    Moreover,
    \[
        \tilde\pi^k_{T(x)}
        =T_\#\Bigl(\frac{h(y)}{h(x)}\pi^k_x(\rmd y)\Bigr),
        \qquad
        \int z\,\tilde\pi^k_{T(x)}(\rmd z)=T(x),
    \]
    so $\tilde\pi^k\in\Mcal(\tilde\mu,\tilde\nu)$.
    Both marginals are supported on $K:=\{(s,t)\in[0,\infty)^2: s+t\le1\}$.
    Denoting $\{z_i\}_i$ the atoms of $\mu$, then $\{T(z_i)\}_i$ are the atoms of $\tilde\mu$, so $\tilde\mu$ is countably supported.
    Moreover, $\tilde\nu$ has positive density on the open sets $T(C)$ and hence full affine span.

    By the proof of \Cref{lem:irreducible.jensen-gaps}, there exists $\chi\in\Mcal(\mu,\nu)$ with $\chi_x\sim\nu$ at every atom $x$ of $\mu$.
    Applying the same transformation to $\chi$ gives $\tilde\chi\in\Mcal(\tilde\mu,\tilde\nu)$ with
    \[
        \tilde\chi_{T(x)} =T_\#\Bigl(\frac{h(y)}{h(x)}\chi_x(\rmd y)\Bigr)
    \]
    Thus $(\tilde\mu,\tilde\nu)$ is irreducible.
    Define $\tilde f$ on $K$ by
    \[
        \tilde f(s,t):=\sup_{j\in\mathbb Z}\max\Bigl\{
        2^js-4^j(1-s-t),\,
        2^jt-2\cdot4^j(1-s-t)
        \Bigr\},
    \]
    and set $\tilde f=+\infty$ outside $K$.
    Then $\tilde f$ is nonnegative, convex and lower semicontinuous, with $\tilde f(T(y))=f(y)/h(y)$.
    In particular, $\tilde f$ is finite at every atom of $\tilde\mu$.
    Furthermore,
    \[
        \int\tilde f\,\rmd\tilde\pi^k_{T(x)}-\tilde f(T(x))
        =\frac{\int f\,\rmd\pi^k_x-f(x)}{h(x)}.
    \]
    Hence $D_{\tilde\pi^k}(\tilde f)=\tfrac1cD_{\pi^k}(f)$ and, by the same argument, $D_{\tilde\chi}(\tilde f)=\tfrac1cD_\chi(f)<\infty$.
    Consequently,
    \[
        D_{\tilde\pi^1}(\tilde f)>D_{\tilde\pi^2}(\tilde f)=0,
        \qquad
        \int\tilde f\,\rmd\tilde\mu
        =\frac1c\int f\,\rmd\mu=+\infty,
        \qquad
        \int\tilde f\,\rmd\tilde\nu
        =\frac1c\int f\,\rmd\nu=+\infty. 
        \qedhere
    \]
\end{proof}

Finally, the compact planar example extends to higher dimensions by taking products with an auxiliary Bass coupling.

\begin{proof}[Proof of \Cref{lem:jensen-gaps}]
    The case $d=2$ follows from \Cref{lem:compact.irreducible.jensen-gaps}.
    For $d>2$, let $\mu,\nu,\pi^1,\pi^2,f$ be as in that lemma.
    Let $\zeta$ be the Bass coupling with $B_0$ uniform on $\{-1,1\}^{d-2}$ and potential $u(b)=\sum_{i=1}^{d-2}|b_i|$.
    Its initial and terminal laws $\beta$ and $\eta$ are uniform on $\{-c,c\}^{d-2}$ and $\{-1,1\}^{d-2}$, respectively, where $c=2\Phi(1)-1$ and $\Phi$ is the standard Gaussian distribution function.
    Moreover, $\zeta_z\sim\eta$ for every atom $z$ of $\beta$.

    Set $\tilde\mu:=\mu\otimes\beta$, $\tilde\nu:=\nu\otimes\eta$, $\tilde\pi^k:=\pi^k\otimes\zeta$ for $k \in \{1,2\}$, and define $\tilde f(x,z):=f(x)$.
    Then $\tilde\pi^k\in\Mcal(\tilde\mu,\tilde\nu)$, $D_{\tilde\pi^k}(\tilde f)=D_{\pi^k}(f)$, and both marginal integrals of $\tilde f$ are infinite.
    By the proof of \Cref{lem:compact.irreducible.jensen-gaps}, there exists $\chi\in\Mcal(\mu,\nu)$ with $\chi_x\sim\nu$ at every atom of $\mu$ and $D_\chi(f)<\infty$.
    Hence $(\chi\otimes\zeta)_{(x,z)}=\chi_x\otimes\zeta_z\sim\tilde\nu$, which proves irreducibility.
    Moreover, $D_{\chi\otimes\zeta}(\tilde f)=D_\chi(f)<\infty$.
    All remaining assertions follow directly from the construction.
\end{proof}

\section{Not every Bass martingale is stretched Brownian motion}\label{sec:Bass!=sBM}

We first prove \Cref{thm:main} in dimension two.
For the next three lemmas, fix $\mu,\nu,\pi^1,\pi^2,f$ as in \Cref{lem:compact.irreducible.jensen-gaps}, and write $\pi^{\rm sBM}:=\operatorname{sBM}(\mu,\nu)$.
The idea is to penalize the Jensen gap in the martingale Benamou--Brenier problem.
For a sufficiently large penalty, the optimizer differs from $\pi^{\rm sBM}$, while dual attainment shows that it nevertheless has Bass structure.
For $\lambda\ge0$, consider the weak martingale transport cost
\[
    C_\lambda(x,\rho)
    :=
    \begin{cases}
        \lambda\bigl(\int f\,\rmd\rho-f(x)\bigr)
        -\MCov(\rho,\gamma),
        &\int y\,\rho(\rmd y)=x,\quad f(x)<\infty,\\
        +\infty,&\text{otherwise},
    \end{cases}
\]
and the corresponding optimization problem
\[
    P_\lambda(\mu,\nu) :=\sup_{\pi\in\Mcal(\mu,\nu)} \bigl\{J(\pi)-\lambda D_\pi(f)\bigr\}.
\]
For $\lambda=0$, this reduces to $\mBB(\mu,\nu)$, whose unique optimizer is $\pi^{\rm sBM}$.

\begin{lemma}\label{lem:existence.kappa}
    There exists $\lambda>0$ such that $P_\lambda(\mu,\nu)$ admits a unique optimizer $\kappa\in\Mcal(\mu,\nu)$ with $\kappa\neq\pi^{\rm sBM}$.
\end{lemma}

\begin{proof}
    Let $K:=\{(s,t) \in[0,\infty)^2: s+t\le1\}$ and $M:=\int|z|\,\gamma(\rmd z)$.
    Since $|y\cdot z|\le |y||z|$, $|\MCov(\rho,\gamma)|\le M$ for every $\rho\in\Pcal(K)$.
    Moreover, $C_\lambda$ is convex in its second argument and jointly lower semicontinuous on $\operatorname{int}K\times\Pcal(K)$, since $f$ is continuous on $\operatorname{int}K$ and lower semicontinuous on $K$, while $\MCov(\cdot,\gamma)$ is continuous on $\Pcal(K)$.
    Jensen's inequality gives $C_\lambda(x,\rho)\ge-M$, while $D_\pi(f)\ge0$ and $D_{\pi^2}(f)=0$ yield
    \[
        0\le J(\pi^2)\le P_\lambda(\mu,\nu)\le M<\infty.
    \]
    By \cite[Theorem~2.2(i)]{BePaRiSc25}, if a weak transport cost $C(x,\rho)$ is convex and lower semicontinuous in $\rho$ for the $p$-weak topology, satisfies $C(x,\rho)\ge a(x)+\int b\,\rmd\rho$ for some $a\in L^1(\mu)$ and $b\in L^1(\nu)$, and the corresponding weak transport problem has finite value, then this infimum is attained by a coupling of $\mu$ and $\nu$.
    Thus $P_\lambda(\mu,\nu)$ admits an optimizer $\kappa \in \Mcal(\mu,\nu)$.
    Uniqueness of $\kappa$ follows from strict concavity of $\MCov(\cdot,\gamma)$, a consequence of the uniqueness of the monotone optimal transport map from the Gaussian law.

    We next show that $D_{\pi^{\rm sBM}}(f)>0$.
    Recall from \Cref{lem:compact.irreducible.jensen-gaps} that $x_0:=T(P_1)=(3/8,1/2)$ is an atom of $\mu$.
    The affine functions $2s+t-1$, $6s+4t-4$ and $2s+3t-2$ lie below $f$ and coincide only at $x_0$.
    Thus their average $\ell(s,t):=\tfrac{10}{3}s+\tfrac83t-\tfrac73$ satisfies $f\ge\ell$, with equality only at $x_0$.
    Since $(\mu,\nu)$ is irreducible, $\pi^{\rm sBM}$ admits a Bass representation, and \Cref{prop:Bass.structure} gives $\pi^{\rm sBM}_{x_0}\sim\nu\neq\delta_{x_0}$.
    In particular,
    \[
        D_{\pi^{\rm sBM}}(f)
        \ge
        \mu(\{x_0\})\int(f-\ell)\,\rmd\pi^{\rm sBM}_{x_0}>0.
    \]
    Since $D_{\pi^2}(f)=0$, choosing $\lambda>0$ sufficiently large gives
    \[
        J(\pi^{\rm sBM})-\lambda D_{\pi^{\rm sBM}}(f)
        <J(\pi^2)=J(\pi^2)-\lambda D_{\pi^2}(f).
    \]
    Hence the unique maximizer $\kappa$ of $P_\lambda(\mu,\nu)$ differs from $\pi^{\rm sBM}$.
\end{proof}

It remains to show that the optimizer $\kappa\neq\pi^{\rm sBM}$ is a Bass coupling.
We first derive a dual formula for $P_\lambda(\mu,\nu)$ and show that its value is attained by a convex potential.
The subsequent lemma then uses this potential to establish \eqref{eq:Bass.structure}.

\begin{lemma}[Strong duality]\label{lem:strong.duality}
    Fix $\lambda>0$ and let $\kappa$ denote the unique optimizer of $P_\lambda(\mu,\nu)$.
    Let $K:=\{(s,t)\in[0,\infty)^2:s+t\le1\}$ and denote by $C(K)$ the space of continuous real-valued functions on $K$.
    For $h\in C(K)$, set $w_h:=(\lambda f+h)^{**}$, extending $\lambda f+h$ by $+\infty$ outside $K$.
    Then
    \begin{equation}\label{eq:penalized.dual}
        P_\lambda(\mu,\nu)
        =
        \inf_{h\in C(K)}
        \Bigl\{
            \int(w_h-\lambda f)\,\rmd\nu
            +
            \int(\lambda f-w_h^C)\,\rmd\mu
        \Bigr\}.
    \end{equation}
    Moreover, there exists a proper lower semicontinuous convex function $w$, finite on $\operatorname{int}\operatorname{co} \operatorname{supp} \nu$, such that
    \[
        P_\lambda(\mu,\nu) = \Ecal_\kappa(w) - \lambda D_\kappa(f).
    \]
\end{lemma}

\begin{proof}
    By \cite[Theorem~1.3]{BaBePa18}, a bounded-below, jointly lower semicontinuous weak transport cost $C$, convex in its second argument, satisfies
    \[
        \inf_{\pi\in\Pi(\mu,\nu)} \int C(x,\pi_x)\,\mu(\rmd x)
        =
        \sup_{h\in C(K)} \Bigl\{\int R_Ch\,\rmd\mu-\int h\,\rmd\nu\Bigr\},
    \]
    where $R_Ch(x):=\inf_{\rho\in\Pcal(K)} \{\int h\,\rmd\rho+C(x,\rho)\}$.
    Applying this to $C_\lambda$, whose assumptions were verified in \Cref{lem:existence.kappa}, yields
    \[
        P_\lambda(\mu,\nu)
        =
        \inf_{h\in C(K)} \Bigl\{\int h\,\rmd\nu-\int R_{C_\lambda}h\,\rmd\mu\Bigr\}.
    \]
    Fix $h\in C(K)$ and recall that $w_h=(\lambda f+h)^{**}$.
    By the definition of $C_\lambda$,
    \begin{equation}\label{eq:proof:penalized.dual:c.conjugate}
        R_{C_\lambda}h(x)+\lambda f(x)
        =
        \inf_{\substack{Y\in K,\;Z\sim\gamma\\ \E[Y]=x}}
        \E[(\lambda f+h)(Y)-Y\cdot Z].
    \end{equation}
    By \cite[Proposition~5.2]{BaBeScTs26}, if $w$ is proper, lower semicontinuous and convex, and
    \[
        \inf_{\rho:\int y\,\rmd\rho=x_0}
        \Bigl\{\int w\,\rmd\rho-\MCov(\rho,\gamma)\Bigr\}>-\infty
    \]
    for some $x_0\in\operatorname{int}(\operatorname{dom}w)$, then
    \[
        \inf_{\substack{\rho\in\Pcal(K)\\
        \int y\,\rho(\rmd y)=x}}
        \Bigl\{\int w\,\rmd\rho-\MCov(\rho,\gamma)\Bigr\}
        =
        \sup_{\zeta\in\R^2}
        \Bigl\{\langle\zeta,x\rangle
        -\int w^*(\zeta+z)\,\gamma(\rmd z)\Bigr\}
        =w^C(x)
    \]
    for every $x\in\operatorname{int}(\operatorname{dom}w)$.
    This applies to $w_h$, since, with $M:=\int|z|\,\gamma(\rmd z)$, the infimum is bounded below by $\min_Kh-M$ and above by $w_h(x)<\infty$ for every $x\in\operatorname{int}K$.

    Although $\lambda f+h$ need not be convex, the right-hand side of \eqref{eq:proof:penalized.dual:c.conjugate} is convex in $x$, finite and hence continuous on $\operatorname{int}K$.
    Its convex conjugate is
    \[
        \int(\lambda f+h)^*(\zeta+z)\,\gamma(\rmd z)
        =\int w_h^*(\zeta+z)\,\gamma(\rmd z).
    \]
    Biconjugacy therefore gives $R_{C_\lambda}h=w_h^C-\lambda f$ on $\operatorname{int}K$, and hence
    \[
        P_\lambda(\mu,\nu)
        =
        \inf_{h\in C(K)} \Bigl\{\int h\,\rmd\nu+ \int(\lambda f-w_h^C)\,\rmd\mu\Bigr\}.
    \]
    Moreover,
    \[
        \lambda f+\min_Kh\le w_h\le\lambda f+\max_Kh,
        \qquad 0\le w_h-w_h^C\le M,
    \]
    so $w_h-\lambda f$ and $\lambda f-w_h^C$ are bounded on $\operatorname{int} K$.
    On the other hand, for every $\pi\in\Mcal(\mu,\nu)$ with $D_\pi(f)<\infty$, the definition of $w_h^C$ gives
    \begin{equation}\label{ineq:proof:penalized.dual:Fenchel}
        w_h^C(x)\le\int w_h\,\rmd\pi_x-\MCov(\pi_x,\gamma),
        \qquad \mu\text{-a.e. }x.
    \end{equation}
    Integrating and using $w_h-\lambda f\le h$, we obtain
    \[
        P_\lambda(\mu,\nu)
        \le
        \int(w_h-\lambda f)\,\rmd\nu+
        \int(\lambda f-w_h^C)\,\rmd\mu
        \le
        \int h\,\rmd\nu+
        \int(\lambda f-w_h^C)\,\rmd\mu.
    \]
    Taking the infimum over $h\in C(K)$ proves \eqref{eq:penalized.dual}.
    
    To prove dual attainment, choose a minimizing sequence $(w_n)_n$ in \eqref{eq:penalized.dual} such that
    \[
        d_n:= \Bigl(\int(w_n-\lambda f)\,\rmd\nu +\int(\lambda f-w_n^C)\,\rmd\mu\Bigr) \ \Big\downarrow \ P_\lambda(\mu,\nu).
    \]
    For every $\pi\in\Mcal(\mu,\nu)$ with $D_\pi(f)<\infty$, cancellation of the integrable differences gives
    \[
        \Ecal_\pi(w_n)=d_n+\lambda D_\pi(f).
    \]
    Let $m:=\int x\,\mu(\rmd x)\in\operatorname{int}K$ and choose $a_n\in\partial w_n^C(m)$.
    Then $\ell_n(x):=w_n^C(m)+a_n\cdot(x-m)$ satisfies $w_n^C\ge\ell_n$, with equality at $m$.
    Subtracting $\ell_n$ does not change $\Ecal_\pi(w_n)$, since $(w_n-\ell_n)^C=w_n^C-\ell_n$.
    Thus we may assume
    \begin{equation}\label{eq:proof:penalized.dual:w_n.normalized}
        w_n\ge w_n^C\ge0,\qquad w_n^C(m)=0.
    \end{equation}
    Moreover, as shown in the proof of \Cref{lem:compact.irreducible.jensen-gaps}, there exists $\chi\in\Mcal(\mu,\nu)$ with $\chi_x\sim\nu$ for $\mu$-a.e.~$x$ and $D_\chi(f)<\infty$.
    Consequently,
    \begin{equation}\label{eq:proof:penalized.dual:Ecal.bounded}
        \sup_n\Ecal_\chi(w_n) =\sup_n d_n + \lambda D_\chi(f)<\infty.
    \end{equation}
    By \cite[Lemma~3.9]{HaJoLoObPa26}, for irreducible $(\mu,\nu)$ and a martingale coupling $\chi$ whose conditional laws charge every open set of positive $\nu$-measure, proper lower semicontinuous convex functions $u_n\ge v_n\ge0$ satisfying
    \[
        v_n(m)=0,\qquad
        \operatorname{int}\operatorname{co}(\operatorname{supp}\nu)
        \subseteq\operatorname{dom}u_n,
        \qquad
        \sup_n\int\Bigl(\int u_n\,\rmd\chi_x-v_n(x)\Bigr)\,\mu(\rmd x)<\infty
    \]
    are locally bounded on $\operatorname{int}\operatorname{co}(\operatorname{supp}\nu)$.
    These conditions hold for $u_n=w_n$ and $v_n=w_n^C$ by \eqref{eq:proof:penalized.dual:w_n.normalized} and \eqref{eq:proof:penalized.dual:Ecal.bounded}.
    Moreover, \cite[Theorem~7.6]{RoWe98} shows that every such sequence admits an epi-convergent subsequence. Denote its limit by $w$.
    Since $w_n\ge0$ and $(w_n)_n$ is locally bounded on $\operatorname{int}\operatorname{co}(\operatorname{supp}\nu)$, the function $w$ is proper, lower semicontinuous, convex and finite on this set.
    Finally, the proof of \cite[Proposition~3.12(1)]{HaJoLoObPa26} establishes that, under the above epi-convergence and normalization,
    \[
        \Ecal_\kappa(w)\le\liminf_n\Ecal_\kappa(w_n).
    \]
    Together with \eqref{ineq:proof:penalized.dual:Fenchel} applied to $w$, this gives
    \[
        J(\kappa)\le\Ecal_\kappa(w)
        \le 
        P_\lambda(\mu,\nu)+\lambda D_\kappa(f)
        =
        J(\kappa).
        \qedhere
    \]
\end{proof}

The dual equality in \Cref{lem:strong.duality} implies that $\kappa_x$ minimizes the expression defining $w^C(x)$ for $\mu$-almost every $x$.
This yields the Bass representation.

\begin{lemma}\label{lem:kappa.is.Bass}
    Let $\kappa$ be the optimizer from \Cref{lem:existence.kappa} and $w$ the convex function obtained in \Cref{lem:strong.duality}.
    Then $\kappa$ is a Bass coupling with potential $w^*$.
\end{lemma}

\begin{proof}
    Let $w$ be as in \Cref{lem:strong.duality}.
    Since $\operatorname{dom}w\subseteq K$, the function $w^*$ is finite and globally Lipschitz, so its gradient is bounded.
    As in the proof of \Cref{prop:Bass.structure}, $w^C$ is differentiable on $\operatorname{int}\operatorname{co}(\operatorname{supp}\nu)$.

    By \Cref{lem:strong.duality}, $w$ satisfies $J(\kappa) - \lambda D_\kappa(f) = \Ecal_\kappa(w) - \lambda D_\kappa(f)$.
    Together with \eqref{ineq:proof:penalized.dual:Fenchel}, this yields
    \[
        \int w\,\rmd\kappa_x-w^C(x)
        =
        \MCov(\kappa_x,\gamma),
        \qquad \mu\text{-a.e. }x.
    \]
    Thus $\kappa_x$ minimizes $\int w\,\rmd\rho-\MCov(\rho,\gamma)$ among probabilities with barycenter $x$.
    According to \cite[Proposition~5.3(ii)]{BaBeScTs26}, every such minimizer with finite value and $x\in\operatorname{int}(\operatorname{dom}w)$ is of the form $(\nabla w^*(\cdot+\zeta))_\#\gamma$, where $\zeta$ maximizes $\zeta\cdot x-\int w^*(\zeta+z)\,\gamma(\rmd z)$.
    Since $w^C=(w^*\star\gamma)^*$ is differentiable, $\zeta=\nabla w^C(x)$, and hence $\kappa$ is a Bass coupling.
\end{proof}

It remains to extend the result to arbitrary $d\ge2$.

\begin{proof}[Proof of \Cref{thm:main}]
    The case $d=2$ follows from \Cref{lem:existence.kappa,lem:kappa.is.Bass}.
    For $d>2$, let $\zeta\in\Mcal(\beta,\eta)$ be the auxiliary Bass coupling from the proof of \Cref{lem:jensen-gaps}, and let $\kappa\in\Mcal(\mu,\nu)$ be the planar non-optimal Bass coupling with potential $w^*$ and Bass measure $\alpha$.
    Since $\zeta$ has Bass potential $u(b)=\sum_{i=1}^{d-2}|b_i|$ and Bass measure $\eta$, the product $\tilde\kappa:=\kappa\otimes\zeta$ is a Bass coupling with potential $v(a,b):=w^*(a)+u(b)$ and Bass measure $\alpha\otimes\eta$.
    Its marginals $\tilde\mu=\mu\otimes\beta$ and $\tilde\nu=\nu\otimes\eta$ have full affine span and are irreducible, since $\kappa_x\sim\nu$ by \Cref{prop:Bass.structure} and $\zeta_z\sim\eta$.
    
    Let $\pi^{\rm sBM}=\operatorname{sBM}(\mu,\nu)$ and set $\tilde\pi:=\pi^{\rm sBM}\otimes\zeta\in\Mcal(\tilde\mu,\tilde\nu)$.
    For every $p\in\Pcal_2(\R^2)$ and $q\in\Pcal_2(\R^{d-2})$, maximal covariance is additive:
    \[
        \MCov(p\otimes q,\gamma_d)
        =
        \MCov(p,\gamma_2)+\MCov(q,\gamma_{d-2}).
    \]
    Consequently, $J(\tilde\pi)-J(\tilde\kappa) = J(\pi^{\rm sBM})-J(\kappa)>0$, and so $\tilde\kappa$ is not stretched Brownian motion.
\end{proof}

\section{Not every cyclically monotone coupling has Bass structure}
\label{sec:mc-not-bass}

In the preceding section, we used a Jensen-gap penalty to construct a Bass coupling which is not stretched Brownian motion.
We now use the same examples from \Cref{lem:jensen-gaps}, but maximize $J$ subject to the constraint $D_\pi(f)=0$.
We show that the resulting couplings are MC-cyclically monotone, although they do not have Bass structure~\eqref{eq:Bass.structure}.

\begin{theorem}\label{thm:mc-not-bass}
    For every $d\ge2$, there exist compactly supported probability measures $\mu\le_{\mathrm{cx}}\nu$ on $\R^d$ with full affine span such that $(\mu,\nu)$ is irreducible and admits an MC-cyclically monotone coupling $\kappa\in\Mcal(\mu,\nu)$ which is not a Bass coupling.
\end{theorem}

We first show that maximizing $J$ among martingale couplings with vanishing Jensen gap yields MC-cyclical monotonicity.

\begin{lemma}\label{lem:mc.zero.gap}
    Let $\mu,\nu\in\Pcal(\R^d)$ be compactly supported, with $\mu$ countably supported, and let $g$ be a non-negative and lower semicontinuous convex function on $\Rd$, finite at every atom of $\mu$.
    If there exists $\pi\in\Mcal(\mu,\nu)$ with $D_\pi(g)=0$, then
    \[
        \sup\{J(\pi):\pi\in\Mcal(\mu,\nu),\ D_\pi(g)=0\}
    \]
    is attained, and every optimizer is MC-cyclically monotone.
\end{lemma}

\begin{proof}
    Let $K$ be a compact set containing both supports and define
    \[
        C(x,\rho):=
        \begin{cases}
            -\MCov(\rho,\gamma_d),
            &\int y\,\rho(\rmd y)=x,\quad \int g\,\rmd\rho\le g(x)<\infty,\\
            +\infty,&\text{otherwise},
        \end{cases}
    \]
    for $\rho\in\Pcal(K)$.
    By Jensen's inequality, the condition $\int g\,\rmd\rho\le g(x)$ is an equality whenever $\rho$ has barycenter $x$.
    Moreover, $C$ is measurable, convex and lower semicontinuous in $\rho$, and bounded below by $-\max_{y\in K}|y|\int|z|\,\gamma_d(\rmd z)$.

    By \cite[Theorem~2.2(i)]{BePaRiSc25}, a weak transport problem with such a cost, bounded below by integrable marginal terms and having finite infimum, admits an optimizer.
    Since there exists $\pi$ such that $D_\pi(g)=0$ and $J$ is bounded, this yields a maximizer $\xi\in\Mcal(\mu,\nu)$ with $D_\xi(g)=0$.

    Moreover, \cite[Theorem~5.3]{BaBePa18} establishes that every optimizer with finite value of a weak transport problem with jointly measurable cost, convex and lower semicontinuous in the conditional law, is $C$-monotone.
    Thus, for any finite family of source atoms $x_1,\ldots,x_N$ and probabilities $q_1,\ldots,q_N$ with barycenters $x_i$ and $\sum_iq_i=\sum_i\xi_{x_i}$, we have
    \begin{equation}\label{eq:proof:mc.zero.gap:mc.monotone}
        \sum_{i=1}^N C(x_i,\xi_{x_i}) \le \sum_{i=1}^N C(x_i,q_i)
    \end{equation}
    Since $D_\xi(g)=0$ and $\sum_iq_i=\sum_i\xi_{x_i}$, we have $\sum_i (\int g\,\rmd q_i-g(x_i))=0$.
    Each summand is nonnegative by Jensen's inequality, hence $\int g\,\rmd q_i=g(x_i)$ for every $i$.
    It follows that $C(x_i,q_i)=-\MCov(q_i,\gamma)$, and so \eqref{eq:proof:mc.zero.gap:mc.monotone} is precisely MC-cyclical monotonicity \eqref{eq:mc-monotonicity}.
\end{proof}

\begin{proof}[Proof of \Cref{thm:mc-not-bass}]
    Fix $d\ge2$ and let $\mu,\nu,\pi^2,f$ be as in \Cref{lem:jensen-gaps}.
    Since $D_{\pi^2}(f)=0$, \Cref{lem:mc.zero.gap} yields an MC-cyclically monotone $\xi\in\Mcal(\mu,\nu)$ with $D_\xi(f)=0$.

    Let $x_0:=T(P_1)$, and recall from the proof of \Cref{lem:existence.kappa} that there exists an affine function $\ell$ on $\R^2$ such that $f(u,z)\ge\ell(u)$, with equality only if $u=x_0$.
    Fix any atom $(x_0,z)$ of $\mu$.
    Since $D_\xi(f)=0$, we have
    \[
        0=\int(f(a,b)-\ell(a))\,\xi_{(x_0,z)}(\rmd a,\rmd b),
    \]
    Thus $\xi_{(x_0,z)}$ is concentrated on $\{x_0\}\times\R^{d-2}$. 
    As $\nu$ has full affine span, $\xi$ cannot be a Bass coupling by \Cref{prop:Bass.structure}.
\end{proof}

\section{Instability of martingale transport beyond the line}\label{sec:instability}

We now prove \Cref{thm:instability} by truncating the Bass measure of the non-optimal coupling constructed in the proof of \Cref{thm:main}.

\begin{proof}[Proof of \Cref{thm:instability}]
    Let $\kappa\neq\pi^{\rm sBM}:=\operatorname{sBM}(\mu,\nu)$ be the Bass coupling from \Cref{thm:main}, with potential $w^*$ and Bass measure $\alpha$.
    For sufficiently large $n$, define
    \[
        \alpha_n(\rmd a) := \frac{\mathbf1_{\{|a|\le n\}}}{\alpha(\{|a|\le n\})}\,\alpha(\rmd a).
    \]
    Let $\kappa^n$ be the Bass coupling with potential $w^*$ and Bass measure $\alpha_n$, and denote its marginals by $\mu^n,\nu^n$.
    Since $\alpha_n\to\alpha$ in total variation, the same holds for $\kappa^n\to\kappa$ and both marginals.
    Moreover, since $\mu,\nu$ have full affine span, so do $\mu^n,\nu^n$ for all sufficiently large $n$.

    Because $\alpha_n\ll\alpha$, we have $\mu^n\ll\mu$, and \eqref{eq:Bass.structure} gives $\kappa^n_x=\kappa_x$ for $\mu^n$-a.e.~$x$.
    Moreover, $\kappa^n_x\sim\nu^n$ and $\kappa_x\sim\nu$, so $\nu^n\sim\nu$ and $(\mu^n,\nu^n)$ is irreducible.
    All marginals are supported on a common compact set, so their convergence also holds in $\Wcal_1$.
    Since $\alpha_n$ has bounded support, \Cref{lem:integrable-bass} gives $w\in L^1(\nu^n)$ and $\kappa^n=\operatorname{sBM}(\mu^n,\nu^n)$.

    As $\kappa^n_x=\kappa_x$ for $\mu^n$-a.e.~$x$ and $\MCov(\kappa_x,\gamma)$ is uniformly bounded on the common compact terminal support, total-variation convergence yields
    \[
        J(\kappa^n)\longrightarrow J(\kappa)<J(\pi^{\rm sBM}),
    \]
    where the strict inequality follows from \Cref{thm:main} and uniqueness of the optimizer.
    Therefore $\kappa^n\to\kappa\neq\operatorname{sBM}(\mu,\nu)$, although each $\kappa^n$ is stretched Brownian motion.

    Finally, let $K:=\operatorname{co}(\operatorname{supp}\nu)$ and suppose that there exist $\pi^n\in\Mcal(\mu^n,\nu^n)$ converging weakly to $\pi^{\rm sBM}$.
    By \Cref{prop:Bass.structure}, $\mu^n$ and $\mu$ are concentrated on $\operatorname{int}K$, where $w^C$ is continuously differentiable.
    Since $w\ge w^C$, the function
    \[
        (x,y)\longmapsto w(y)-w^C(x)-\nabla w^C(x)\cdot(y-x),
        \qquad x\in\operatorname{int}K,
    \]
    is nonnegative and lower semicontinuous.
    Its integral against a martingale coupling $\pi$ with first marginal concentrated on $\operatorname{int}K$ equals $\Ecal_\pi(w)$.
    Hence Portmanteau's theorem gives
    \[
        \Ecal_{\pi^{\rm sBM}}(w)
        \le\liminf_{n\to\infty}\Ecal_{\pi^n}(w).
    \]
    Since $w\in L^1(\nu^n)$, \eqref{eq:bass-conditional-equality} also gives $w^C\in L^1(\mu^n)$.
    Thus, for every $\pi^n\in\Mcal(\mu^n,\nu^n)$,
    \[
        \Ecal_{\pi^n}(w)
        =\int w\,\rmd\nu^n-\int w^C\,\rmd\mu^n
        =\Ecal_{\kappa^n}(w)
        =J(\kappa^n).
    \]
    Together with Inequality~\eqref{eq:conditional-dual-bound}, this implies
    \[
        J(\pi^{\rm sBM})
        \le\Ecal_{\pi^{\rm sBM}}(w)
        \le\liminf_{n\to\infty}\Ecal_{\pi^n}(w)
        =J(\kappa)<J(\pi^{\rm sBM}),
    \]
    a contradiction.
\end{proof}

\section*{Acknowledgments}
The authors thank Walter Schachermayer for many stimulating discussions and valuable suggestions during the preparation of this manuscript.

The authors further acknowledge the substantial assistance of OpenAI's GPT-6 Astra in obtaining a preliminary version of the construction in \Cref{sec:construction}, following extensive work under the authors' direction. The authors subsequently verified, substantially revised, and extended the construction, which forms the basis for the results of this paper.

\printbibliography

\end{document}